\documentclass[]{iopart}
\usepackage{iopams}
 \usepackage{amssymb}
\usepackage{setstack}
\usepackage{latexsym}
\usepackage[usenames]{color}

\newtheorem{theorem}{Theorem}[section]
\newtheorem{corollary}[theorem]{Corollary}
\newtheorem{defn}[theorem]{Definition}
\newtheorem{lem}[theorem]{Lemma}
\newtheorem{proposition}[theorem]{Proposition}
\newtheorem{problem}[theorem]{Problem}
\newtheorem{remark}[theorem]{Remark}

\begin{document}
\eqnobysec
\title{KdV equation under periodic boundary conditions and its perturbations}
\author{HUANG Guan$^1$ and Sergei KUKSIN$^2$}
\address{$^1$ C.M.L.S, \'Ecole Polytechnique, Palaiseau, France}
\address{$^2$ I.M.J, Universit\'e Paris Diderot-Paris 7, Paris, France}

\begin{abstract}
In this paper we discuss properties of the KdV equation under  periodic boundary conditions, especially those which are important to study  perturbations of the equation. Next we  review what is known now about long-time behaviour of solutions for  perturbed KdV equations.
\end{abstract}
\tableofcontents

\maketitle
\bibliographystyle{plain}
\setcounter{section}{-1}

\section{Introduction}
The famous  Korteweg-de Vries (KdV) equation \[u_t=-u_{xxx}+6uu_x, \quad x\in\mathbb{R},\]
was first proposed by Joseph Boussinesq \cite{bou1871} as a model for shallow water wave propagation. 
It became famous later when   two Dutch mathematicians, Diederik  Korteweg and  Gustav De~Vries \cite{kdv1895}, used it 
to explain the existence of a soliton water wave, previously observed by John  Russel in physical experiments. Their  work was so successful that this  equation is now named after them. Since the mid-sixties of 20th century the KdV equation received  a lot of attention from mathematical  and physical communities after the numerical results of Kruskal and Zabusky \cite{krz1965}  led to the discovery that  its solitary wave solutions interact in an integrable way. It turns out that in some suitable setting, the KdV equation can be viewed as an integrable infinite dimensional hamiltonian  system.

In his ``New Methods of Celestial Mechanics",  Poincar\'e calls the task to study perturbations of integrable systems the ``General Problem of Dynamics".  The great scientist was motivated by the celestial mechanics, where perturbed integrable systems play a very important role\footnote{For example, the solar system, regarded as a system of 8 interacting planets rotating around the Sun, is a small perturbation of the Kaplerian system. The latter is integrable}. For a similar reason his maxim is true for mathematical physics, where many important processes are described by suitable perturbations of an integrable PDE, while the unperturbed integrable equations correspond to idealization of physical reality. In particular, no physical process is exactly described by the KdV equation. 

In this paper\footnote{ Based on the courses, given by the second author  in Saint-Etienne de Tin\'ee in February  2012 and in the High School for Economics (Moscow) in April 2013.}, we focus on the KdV equation with zero mean value periodic boundary condition. It is known since the works of Novikov, Lax,  Marc\^enko, Its-Matveev and McKean-Trubowitz that this system is integrable (\cite{mar1977, mct1978, lax1968, nov1974}).  All of its solutions are periodic, quasi-periodic or almost periodic in time. In Section~1 we discuss the KdV equation in the framework of infinite-dimensional hamiltonian  systems, in Section 2 we present some normal form results for finite-dimensional integrable hamiltonian  systems and in Section 3 --- recent results on KdV which may be regarded as infinite-dimensional versions of those in Section 2. Finally in Section~4 we discuss long-time behavior of solutions for the perturbed KdV equations, under hamiltonian  and non-hamiltonian  perturbations.  Results presented there  are heavily based on theorems from Section 3.

As indicated  in  the title of our work, we restrict our study to the periodic boundary conditions. In this case the KdV equation behaves as a hamiltonian  system with countablely-many degrees of freedom and the method of Dynamical Systems may be used for its study (same is true  for other hamiltonian  PDEs in finite volume, e.g. see \cite{k06_handb}).
The KdV is a good example of an integrable PDE in the sense that properties of many other integrable equations with self-adjoint Lax operators,  e.g. of the defocusing Zakharov-Shabat equation (see \cite{GK}), 
 and of their perturbations are very similar to those of KdV and its perturbations, while the equations with non-selfadjoint Lax operators, e.g. the Sine-Gordon 
equation, 
 are similar to KdV when we study their small-amplitude solutions (and the KAM-theory for such equation is similar to the KAM theory for  KdV  without the smallness assumption, see \cite{kuk2000}).
When considered on the whole line with ``zero at infinity" boundary condition, due to the effect of radiation, the KdV equation and its perturbations behave differently, and people working on these problems prefer to call them ``dispersive systems". To discuss  the corresponding results should be a topic of another work. 

The number of publications, dedicated to KdV and its perturbations is immense, and our bibliography is hopelessly incomplete.

 \section{KdV under periodic boundary conditions  as a hamiltonian  system}
Consider the KdV equation under zero mean value periodic boundary condition:
\begin{equation}
u_t+u_{xxx}-6uu_x=0,\quad x\in\mathbb{T}=\mathbb{R}/\mathbb{Z},\quad \int_{\mathbb{T}}udx=0.\label{kdv1}
\end{equation}
(Note that the mean-value $\int_{\mathbb{T}}udx$ of  a space-periodic solution $u$ is a time-independent quantity, to simplify presentation we choose it to be zero.) 
To fix the setup,  for any integer $p\geqslant0$, we introduce  the Sobolev space of real valued functions on $\mathbb{T}$ with zero mean-value:
\[
\fl
H^p=\Big\{u\in L^2(\mathbb{T},\mathbb{R}):\;||u||_p<+\infty,\;\int_{\mathbb{T}} =0\Big\},\quad
||u||_p^2=\sum_{k\in\mathbb{N}}|2\pi k|^{2p}(|\hat{u}_k|^2+|\hat{u}_{-k}|^2).\]
Here $\hat{u}_k$, $\hat{u}_{-k}$,  $k\in\mathbb{N}$, are the Fourier coefficients of $u$ with respect to the trigonometric base
\begin{equation}
e_k=\sqrt{2}\cos{2\pi kx},\quad k>0\quad \mbox{and}\quad e_k=\sqrt{2}\sin{2\pi kx},\quad k<0, 
\label{basis1}
\end{equation}
 i.e. 
 \begin{equation}u=\sum_{k\in\mathbb{N}}\hat{u}_ke_k+\hat{u}_{-k}e_{-k}.
 \label{basis2}
 \end{equation}
In particular, $H^0$ is the space of $L^2$-functions on $\mathbb{T}$ with zero mean-value. By $\langle \cdot,\cdot\rangle$ we denote the scalar product in $H^0$ (i.e. the $L^2$-scalar product). 

For a $C^1$-smooth functional $F$ on some space $H^p$, we denote by $\nabla F$ its gradient with respect to $\langle \cdot,\cdot\rangle$, i.e.
\[d F(u)(v)=\langle\nabla F(u),v\rangle,\]
if $u$ and $v$ are sufficiently smooth. So $\nabla F(u)=\frac{\delta F}{\delta u(x)}+const$, where $\frac{\delta F}{\delta u}$ is the variational derivative, and the constant is chosen in such a way that the mean-value of the r.h.s vanishes.  See \cite{kuk2000,kjp2003}  for details.
The initial value problem for  KdV on the circle $\mathbb{T}$ is well posed on every Sobolev space $ H^p$ with $p\geqslant 1$, see \cite{sat1976}.  The regularity of KdV in function spaces of lower smoothness was studied intensively, see  \cite{Tao, KT} and references in these works; also see \cite{Tao} for some qualitative results concerning the KdV flow in these spaces.
We avoid this topic.

 It was observed by Gardner \cite{ gar1971} that if we introduce the Poisson bracket which assigns to any two functionals $F(u)$ and $G(u)$ the new functional $\{F,G\}$,
\begin{equation}
\Big\{ F,G\Big\}(u)=\int_{\mathbb{T}}\frac{d}{dx}\nabla F(u(x))\nabla G(u(x))dx \label{poisson1}
\end{equation}
(we assume that the r.h.s is well defined, see \cite{kuk2000,  k06_handb,  kjp2003} for details), then KdV becomes a hamiltonian  PDE.
Indeed,  this  bracket corresponds  to a differentiable hamiltonian  function $F$  a vector filed $\mathcal{V}_F$, such that 
\[\langle\mathcal{ V}_F(u),\nabla G(u)\rangle=\{F,G\}(u)
\]
for any differentiable functional $G$. From this  relation we see that 
$
\mathcal{V}_F(u)=\frac{\partial }{\partial x}\nabla F(u).
$
So the KdV equation takes the hamiltonian  form
\begin{equation}u_t=\frac{\partial }{\partial x}\nabla \mathcal{H}(u),
\label{okdveh1}
\end{equation}
with the KdV Hamiltonian
\begin{equation}
\mathcal{H}(u)=\int_{\mathbb{T}}(\frac{u_x^2}{2}+u^3)dx.\label{kdvh}
\end{equation}
The Gardner bracket (\ref{poisson1}) corresponds to  the symplectic structure, defined in $H^0$ (as well as in any space $H^p$, $p\geqslant0$) by the 2-form
\begin{equation}\omega_2^{G}(\xi,\eta)=\big\langle(-\frac{\partial}{\partial x})^{-1}\xi,\eta\big\rangle\quad\mbox{for}\quad\xi,\eta\in H^0.
\label{omega1}
\end{equation}
Indeed, since $\omega_2^G(\mathcal{V}_F(u),\xi)\equiv-\langle\nabla F(u),\xi\rangle$, then the 2-form $\omega_2^G$ also assigns to a Hamiltonian  $F$   the vector field $\mathcal{V}_F$ (see \cite{arn1989, kjp2003,kuk2000,k06_handb}).

We note that the bracket (\ref{poisson1}) is well defined on the whole Sobolev spaces $H^p(\mathbb{T})=H^p\oplus\mathbb{R}$, while the symplectic form $\omega_2^G$ is not, and the affine subspaces $\{u\in H^p(\mathbb{T}):\;\int_{\mathbb{T}}udx=const\}\simeq H^p$ are symplectic leaves for this Poisson system. We study the equation only on the leaf $\int_{\mathbb{T}}udx=0$, but on other leaves it may be studied similarly.  

Writing a function $u(x)\in H^0$ as in (\ref{basis2}) we see that 
$\omega_2^G=\sum_{k=1}^{\infty}k^{-1}d\hat{u}_k\wedge\hat{u}_{-k}$
and that 
$\mathcal{H}(u)=H(\hat{u}):=\Lambda(\hat{u})+G(\hat{u})$ 
with 
\[\fl
\qquad\quad
\Lambda(\hat{u})=\sum_{k=1}^{+\infty}(2\pi k)^2\big(\frac{1}{2}\hat{u}_k^2+\frac{1}{2}\hat{u}_{-k}^2\big), 
\qquad G(\hat{u})=\sum_{k,l,m\neq 0,k+l+m=0}\hat{u}_k\hat{u}_l\hat{u}_m.\]
Accordingly, the KdV equation  may be written as the infinite chain of hamiltonian  equations
\[\frac{d}{dt}\hat{u}_j=-2\pi j\frac{\partial H(\hat{u})}{\partial \hat{u}_{-j}},\quad j=\pm1,\;\pm2,\dots.\]

\section{Finite dimensional integrable systems}
Classically, integrable systems are particular hamiltonian  systems that can be integrated in quadratures. It was observed by Liouville that for a hamiltonian  system with $n$ degrees of freedom to be integrable, it has to possess $n$ independent integrals in involution.  This assertion can be understood globally and locally. Now we recall corresponding finite-dimensional definitions and results.
\subsection{Liouville-integrable systems}
Let $Q\subset \mathbb{R}^{2n}_{(p,q)}$ be a $2n$-dimensional domain. We provide it with the standard symplectic form 
$\omega_0=dp\wedge dq$ 
and the corresponding Poisson bracket 
\[\{f,g\}=\nabla_p f\cdot\nabla_qg-\nabla_qf\cdot\nabla_pg,\]
where $g,f\in C^1(Q)$ and $``\cdot"$ stands for the Euclidean  scalar product in $\mathbb{R}^n$ (see \cite{arn1989}).  If $\{f,g\}=0$, the functions $f$ and $g$ are called {\it commuting}, or {\it in involution}. If $h(p,q)$ is a $C^1$-function on $Q$, then the hamiltonian  system with the Hamiltonian $h$ is 
\begin{equation}
\dot{p}=-\nabla_q h,\quad \dot{q}=\nabla_p h.\label{hfinite}
\end{equation}
\begin{defn} (Liouville-integrability). The hamiltonian  system (\ref{hfinite}) is called integrable in the sense of Liouville if  its Hamiltonian $h$ admits $n$ independent integrals in involution $h_1,\dots,h_n$. That is, 
 $\ \{h,h_i\}=0$ for $1\leqslant i\leqslant n$; 
 $\{h_i,h_j\}=0$ for $1\leqslant i,\;j\leqslant n$, and 
  $dh_1\wedge\cdots\wedge dh_n\neq0$.

\end{defn}

A nice structure of an Liouville-integrable system is assured  by the celebrated Liouville-Arnold-Jost theorem 
(see \cite{arn1989, MoS})
It claims that if  an integrable systems  is such that the level sets $T_c=\{(p,q)\in Q: h_1(p,q)=c_1,\dots,h_n(p,q)=c_n\}$, $c=(c_1,\dots,c_n)\in\mathbb{R}^n$ are compact, then each non-empty set $T_c$ is an embedded  $n$-dimensional  torus.  Moreover for a suitable neighborhood $O_{T_c}$ of $T_c$  in $Q$ there exists a symplectomorphism 
\[\Theta:\;O_{T_c}\to O\times\mathbb{T}^n=\{(I,\varphi)\},\quad O\subset\mathbb{R}^n,\]
where the symplectic structure in $O\times\mathbb{T}^n$ is given by the 2-form $dI\wedge d\varphi$. Finally, there exists a function $\bar{h}(I)$ such that 
$ h(p,q)=\bar{h}(\Theta(p,q)).$
This result is true both in the smooth and analytic categories.

The coordinates $(I,\varphi)$ are called the {\it action-angle variables} for (\ref{hfinite}). 
Using them   the hamiltonian  system may be written as 
\begin{equation}\dot{I}=0,\quad\dot{\varphi}=\nabla_I \bar{h}(I).
\label{laj1}
\end{equation}
Accordingly, in the original 
  coordinates $(p,q)$  solutions of the  system  are 
\[(p,q)(t)=\Theta^{-1}(I_0,\varphi_0+\nabla_I\bar{h}(I_0)t).\]

On $O\times\mathbb{T}^n$, consider  the 1-form \mbox{$Id\varphi=\sum_{j=1}^nI_jd\varphi_j$}, then \mbox{$d(Id\varphi)=dI\wedge d\varphi$}. For any vector $I\in O$, and for $j=1,\dots,n$,   denote by  $C_j(I)$ the cycle \mbox{$\{(I,\varphi)\in O\times\mathbb{T}^n:\varphi_j\in[0,2\pi]\;\mbox{and}\;\varphi_i=const,\;\mbox{if}\;i\neq j\}$}.  Then 
\[\frac{1}{2\pi}\int_{C_j}Id\varphi=\frac{1}{2\pi}\int_{C_j}I_jd\varphi_j=I_j.\]
Consider a disc $D_j\subset O\times\mathbb{T}^n$ such that $\partial D_j=C_j$. For any 1-form $\omega_1$,
satisfying 
 $d\omega_1=dI\wedge d\varphi$, we have
\[\frac{1}{2\pi}\int_{C_j}(Id\varphi-\omega_1)=\frac{1}{2\pi}\int_{D_j}d(Id\varphi-\omega_1)=0.\]
So
\begin{equation}
I_j=\frac{1}{2\pi}\int_{C_j(I)}\omega_1,\quad \mbox{if}\quad d\omega_1=dI\wedge d\varphi.\label{arnold1}
\end{equation}
This is  the {\it Arnold formula for actions}.

\subsection{Birkhoff Integrable systems}
We denote by $\mathbb{J}$ the standard symplectic matrix
 $
 \mathbb{J}=\mbox{diag}\Big\{\left(\begin{array}{cc}0 & -1 \\1 & 0\end{array}\right)\Big\},
 $
 operating in any $\mathbb{R}^{2n}$ (e.g. in  $\mathbb{R}^2$).
Assume that the origin is an elliptic critical point of a smooth Hamiltonian $h$, i.e. $\nabla h(0)=0$ and that  the  matrix $\mathbb{J}\nabla^2h(0)$ has only pure imagine eigenvalues. Then there exists a linear symplectic change of coordinates  which puts   $h$ to  the form
\[h=\sum_{i=1}^n\lambda_i(p_i^2+q_i^2)+h.o.t, \qquad \lambda_j\in\mathbb{R} \quad\forall j.\]
 If the frequencies $(\lambda_1,\dots,\lambda_n)$ satisfy some non-resonance conditions, then this normalization process can be carried out to higher order terms.  The result of this normalization is known as the {\it Birkhoff normal form for the Hamiltonian $h$.}
 
\begin{defn}
The frequencies $\lambda_1,\dots,\lambda_n$ are non-resonant up to order $m\geqslant1$ if 
$\sum_{i=1}^nk_i\lambda_i\neq0$ for each $k\in\mathbb Z^n$ such that   $1\leqslant \sum_{i=1}^n|k_i|\leqslant m$. 
 They are called  non-resonant if $k_1\lambda_1+\dots+k_n\lambda_n=0$ with integers $k_1,\dots,k_n$ only when all $k_j$'s 
  vanish. 
\end{defn}
\begin{theorem}
(Birkhoff normal form, see \cite{mos1968, MoS}) Let $H=N_2+\cdots$ be a real analytic Hamiltonian in the vicinity of  the origin in $(\mathbb{R}^{2n}_{(p,q)},dp\wedge dq )$ with the quadratic part $N_2=\sum_{i=1}^n\lambda_i(q_i^2+p_i^2)$. If the (real) frequencies $\lambda_1,\dots,\lambda_n$ are non-resonant up to order~\mbox{$m\geqslant3$}, then there exists a real analytic symplectic trasformation $\Psi_m=Id+\cdots$, such that
\[H\circ\Psi_m=N_2+N_4+\dots+N_m+h.o.t.\]
Here $N_i$ are homogeneous polynomials of order $i$, which are actually smooth functions of  variables $p_1^2+q_1^2,\dots,p_n^2+q_n^2$.  If the frequencies are non-resonant, then there exists a formal symplectic transformation $\Psi=Id+\cdots$, represented by a formal power series, such that  
$H\circ\Psi=N_2+N_4+\cdots$
(this equality holds in the sense of  formal series).
\label{Birkhoff}
\end{theorem}

If the transformation, converting   $H$ to the  Birkhoff normal form,  was convergent, then the resulting Hamiltonian would be integrable in a neighborhood of the origin with the integrals  $p_1^2+q_1^2,\dots,p_n^2+q_n^2$. These functions are not independent when  $p_i=q_i=0$ for some $i$. So the system is not integrable   in the sense of Liouville. But it is integrable in a weaker sense:
\begin{defn}
Functions $f_1,\dots,f_k$ are  functionally independent if their differentials $df_1,\dots,df_k$ are linearly independent on a dense open set.  A $2n$-dimensional Hamiltonian is called Birkhoff integrable near
an equilibrium~\mbox{$m\in\mathbb{R}^{2n}$}, if it admits $n$ functionally independent integrals in involution in the vicinity of $m$.
\end{defn}

Birkhoff normal form provides a powerful tool to  study the dynamics of hamiltonian  PDEs, e.g. see \cite{ kup1996, bmg2006} and references in \cite{ bmg2006} . However, in this paper we shall not discuss its version  for KdV, since for that equation  there exists a stronger normal form.  Now we  pass to its  counterpart in finite dimension. 

\subsection{Vey theorem}
The results of this section hold both in the $C^{\infty}$-smooth and analytic categories.

\begin{defn}
Consider a 
 Birkhoff integrable system, defined 
  near an equilibrium~\mbox{$m\in\mathbb{R}^{2n}$}, with
  independent commuting  integrals $F=(F_1,\dots,F_n)$. Its Poisson algebra
   is the linear space 
$\,
\mathcal{A}(F)=\Big\{G: \{G,F_i\}=0,\;i=1,\dots,n\Big\}.$
\end{defn}
Note that although the integrals of an integrable system 
 are not defined in a unique way, the corresponding algebra $\mathcal{A}(F)$ 
 is.
\begin{defn} A Poisson algebra $\mathcal{A}(F)$  is said to be non-resonant at a point 
$m\in {\mathbb R}^{2n}$, if it contains a Hamiltonian with a non-resonant  elliptic critical point at~$m$ (i.e., around $m$ one can introduce symplectic coordinates $(p,q)$ such that the quadratic part of that Hamiltonian at $m$ is $\sum\lambda_j(p_j^2+q_j^2)$, where the real numbers $\lambda_j$ are non-resonant).
\end{defn}

It is easy to verify that if some $F_1\in\mathcal{A(F)}$   is elliptic and non-resonant 
at the equilibrium~$m$, then all other functions in $\mathcal{A(F)}$ are elliptic at~$m$ as well.  
 
\begin{theorem} (Vey's theorem).
 Let $F=(F_1,\dots,F_n)$ be n functionally  independent functions in involution in a neighbourhood of a point $m\in {\mathbb R}^{2n}$. If the Poisson algebra $\mathcal{A}(F)$ is non-resonant at $m$, then one can introduce  around $m$ symplectic coordinates $(p,q)$ so that $\mathcal{A}(F)$ consists of all functions, which are actually functions of~\mbox{$p_1^2+q_1^2,\dots,p_n^2+q_n^2$}.
\label{th-vey}
\end{theorem}

\noindent{\bf Example.} Let $F=(f_1,\dots, f_n)$ be a system of smooth 
commuting  Hamiltonians, defined in the vicinity of their joint equilibrium 
$m\in\mathbb{R}^{2n}$, such that 
 the hessians $\nabla^2f_i(m)$, $1\leqslant i\leqslant n$, are linear independent. Then the theorem above applies to the 
Poisson algebra  $\mathcal{A}(F)$.

In \cite{vey1978}  Vey proved   the theorem in the analytic case with an additional non-degeneracy condition, which  was later removed by Ito in \cite{ito1989}. The  results in \cite{vey1978,ito1989} also apply to  non-elliptic cases.  The smooth version of Theorem \ref{th-vey} is due to Eliasson \cite{eli1990}. 
There exists an infinite dimensional extension of the theorem, see \cite{kpe2010}.

\section{Integrability of KdV} 
The KdV equation (\ref{kdv1}) admits infinitely many  integrals in involution, and there are different ways to obtain  them, see \cite{gar1971,miu1968,mgk1968, lax1968, zmn1984}. Below we present  an elegant way to  construct a set of Poisson commuting integrals by considering the spectrum of an associated Schr\"odinger operator, due to Piter~Lax  \cite{lax1968} (see \cite{lax1996} for a nice presentation of the theory).

\subsection{Lax pair}\label{s_Lax}
Let $u(x)$ be a  $L^2$-function on $\mathbb{T}$.  Consider the differential operators $L_u$ and  $B_u$, acting on $2$-periodic functions\footnote{note the doubling of the period.}
\[L_u=-\frac{d^2}{dx^2}+u(x),\quad B_u=-4\frac{d^3}{dx^3}+3u(x)\frac{d}{dx}+3\frac{d}{dx}u(x),\]
where we view $u(x)$ as a multiplication operator $f\mapsto u(x)f$. The operators $B_u$ and $L_u$ are called the {\it Lax pair} for KdV. 
  Calculating the commutator $[B_u,L_u]=B_uL_u-L_uB_u$, we  see that most of the terms cancel  and the only term left is $-u_{xxx}+6uu_x$.  Therefore  if $u(t,x)$ is a solution of (\ref{kdv1}), then the operators $L(t)=L_{u(t,\cdot)}$ and $B(t)=B_{u(t,\cdot)}$ satisfy the operator equation
\begin{equation}
\frac{d}{dt}L(t)=[B(t),L(t)].\label{laxpair1}
\end{equation}
Note that the operator $B(t)$ are skew-symmetric, $B(t)^*=-B(t)$. 
Let $U(t)$ be the one-parameter family of unitary operators,  defined by the  differential equation
\[\frac{d}{dt}U=B(t)U, \quad U(0)={\rm Id}.\]
Then 
$\  L(t)=U^{-1}(t)L(0)U(t)$. 
Therefore, the operator  $L(t)$  is unitary conjugated to $L(0)$. Consequently, its spectrum is independent of $t$. That is, the  spectral  data of  the operator $L_u$ provide a set of conserved quantities for the KdV equation (\ref{kdv1}). Since $L_u$ is the strurm-Liouville operator with a potential $u(x)$, then in the context of this theory functions $u(x)$ are called {\it potentials}.

It is well known that for any $L^2$-potential $u$ the spectrum of the Sturm-Liouville operator $L_u$,  regarded as an unbounded operator in $L^2(\mathbb{R}/2\mathbb{Z})$,  is a sequence of simple or double eigenvalues $\{\lambda_j: j\geqslant 0\}$, tending to infinity:
\[\mbox{spec}(u)=\{\lambda_0<\lambda_1\leqslant\lambda_2<\cdots\nearrow \infty\}.
\]
Equality or inequality may occur in every place with a "$\leqslant$" sign (see \cite{mar1977,kjp2003}).  The segment $[\lambda_{2j-1},\lambda_{2j}]$ is called the {\it $n$-th spectral gap}. The asymptotic behaviour of the periodic eigenvalues is 
\[\lambda_{2n-1}(u),\;\lambda_{2n} (u)=n^2\pi^2+[u]+l^2(n),\]
where $[u]$ is the mean value of $u$, and $l^2(n)$ is the $n$-th number of an $l^2$ sequence.
Let $g_n(u)=\lambda_{2n}(u)-\lambda_{2n-1}(u)\geqslant0$, $n\geqslant 1$. These quantities are  conserved  under the flow of KdV. We call $g_n$ the {\it $n$-th gap-length} of the spectrum. The $n$-th gap is called {\it open} if $g_n>0$, otherwise it is  {\it closed}.  However, from  the analytic point of view the periodic eigenvalues and the gap-lengths are not satisfactory integrals, since $\lambda_n$ is not  a smooth function of  the potential $u$ when  $g_n=0$. Fortunately, the squared gap lengths $g_n^2(u)$, $n\geqslant 1$, are real analytic functions on  $L^2$, which  Poisson commute with each other (see \cite{mct1976, lax1996, kjp2003}). 
Moreover, together with the mean value, the  gap lengths  determine uniquely the periodic spectrum of a potential,  and their asymptotic behavior characterizes the regularity of a potential in exactly  the same way as its Fourier coefficients \cite{mar1977, gat1984}.
\medskip

This method applies  to integrate other
  hamiltonian systems in  finite or infinite dimension. It is remarkably general and is referred to as the 
 {\it method of Lax pair.}

\subsection{Action-angle coordinates}
We denote by $\mbox{Iso}(u_0)$   the isospectral set of a potential $u_0\in H^0$:
\[\mbox{Iso}(u_0)=\Big\{u\in H^0: \quad \mbox{spec}(u)=\mbox{spec}(u_0)\Big\}.\]
It is  invariant under the flow of KdV and may   be characterized by the gap lengths
\[\mbox{Iso}(u_0)=\Big\{u\in H^0:\quad g_n(u)=g_n(u_0),\;n\geqslant 1\Big\}.\]
Moreover,  for any $n\geqslant1$, $u_0\in H^n$ if and only if $\mbox{Iso}(u_0)\subset H^n$.

In \cite{mct1976}, McKean and Trobwitz showed that  the $\mbox{Iso}(u_0)$ is   
 homemorphic to a compact  torus, whose dimension equals the number of open gaps. 
  So the phase space $H^0$ is foliated by a collection of KdV-invariant  tori of different dimensions, finite or infinite.   A potential $u\in H^0$ is called {\it finite-gap}  if only a finite number of its spectral gaps are open. 
The finite-dimensional KdV-invariant torus Iso$(u_0)$ is  called a {\it finite-gap torus}. 
For any $n\in\mathbb{N}$ let us  set
\begin{equation}\mathcal{J}^n=\Big\{u\in H^0:\;\;g_j(u)=0\;\;\mbox{if}\;\;j>n\Big\}.
\label{finite-gap}
\end{equation}
 We call the sets  $\mathcal{J}^n$, $n\in\mathbb{N}$, the finite-gap manifolds.
 
 \begin{theorem} For any $n\in\mathbb{N}$, the finite gap manifold
  $(\mathcal{J}^n,\omega^G_2)$  is a smooth symplectic $2n$-manifold, invariant under the flow of KdV (\ref{kdv1}), and 
 \[T_0\mathcal{J}^n=\Big\{u\in H^0:\;\;\hat{u}_k=0\;\;\mbox{if}\;\; |k|\geqslant n+1\Big\},\]
 (see (\ref{basis2})). Moreover, the square gap lengths $g_k^2(u)$, $k=1,\dots,n$, form $n$ commuting analytic integrals of motions, non-degenerated everywhere on the dense domain 
  \mbox{$\mathcal{J}_0^n=\{u\in\mathcal{J}^n:g_1(u),\dots,g_n(u)>0\}$}.
 \end{theorem}
Therefore, the Liouville-Arnold-Jost theorem applies everywhere on $\mathcal{J}_0^n$, $n\in\mathbb{N}$. Furthermore, the union of the 
 finite gap manifolds $\cup_{n\in\mathbb{N}}\mathcal{J}^n$ is dense in each space
 $H^s$  (see \cite{mar1977}). 
This hints that on the spaces $H^s$, $s\geqslant0$, it may be possible to construct global action-angle coordinates   for KdV. In \cite{flm1976}, 
 Flaschka and McLaughlin used the Arnold formula (\ref{arnold1}) to get an  explicit formula for  action variables  of  KdV  in terms of the 2-period spectral data of $L_u$.   To explain their construction, denote by
 $y_1(x,\lambda,u)$ and $y_2(x,\lambda,u)$  the standard fundamental solutions of  the equation
$\ -y^{\prime\prime}+uy=\lambda y$, 
defined by the initial conditions
\begin{eqnarray*}y_1(0,\lambda,u)=1,\quad y_2(0,\lambda,u)=0,\\
y_1^{\prime}(0,\lambda,u)=0,\quad y_2^{\prime}(0,\lambda,u)=1.
\end{eqnarray*}
The quantity $\triangle(\lambda,u)=y_1(1,\lambda,u)+y_2^{\prime}(1,\lambda,u)$ is called the {\it discriminant}, associated with this pair of solutions. The periodic  spectrum of $u$ is precisely the zero set of the entire function $\triangle^2(\lambda,u)-4$, for which we have the explicit representation (see e.g. \cite{zmn1984, mct1976})
\[\triangle^2(\lambda,u)-4=4(\lambda_0-\lambda)\prod_{n\geqslant1}\frac{(\lambda_{2n}-\lambda)(\lambda_{2n-1}-\lambda)}{n^4\pi^4}.
\] 
This function   is  a spectral invariant.  We also need the spectrum of the differential operator $L_u=-\frac{d^2}{dx^2}+u$ under Dirichlet boundary conditions on the interval $[0,1]$. It consists of an unbounded sequence of single Dirichlet eigenvalues
\[\mu_1(u)<\mu_2(u)<\dots\nearrow \infty,\]
which satisfy $\lambda_{2n-1}(u)\leqslant \mu_n(u)\leqslant \lambda_{2n}(u)$, for all $n\in\mathbb{N}$.  Thus, the $n$-th Dirichlet eigenvalue $\mu_n$ is always contained in the $n$-th spectral gap. 
 The Dirichlet spectrum provides  coordinates on the isospectral sets (see \cite{mct1976,mar1977, kjp2003}). For any $z\in\mathbb{T}$, denote by $\{\mu_j(u,z),\;j\geqslant1\}$ the spectrum of the operator $L_u$ under the shifted Dirichlet boundary conditions $y(z)=y(z+1)=0$ (so $\mu_j(u,0)=\mu_j(u)$); still $\lambda_{2n-1}\leqslant \mu_n(u,z)\leqslant\lambda_{2n}(u)$. Jointly with the spectrum $\{\lambda_j\}$, it defines the potential $u(x)$ via the remarkable {\it trace formula} (see \cite{zmn1984, DMN1976,kjp2003,mct1976}):
\[u(z)=\lambda_0(u)+\sum_{j=1}^{\infty}(\lambda_{2j-1}(u)+\lambda_{2j}(u)-2\mu_j(u,z)).\]

Define 
\[f_n(u)=2\log(-1)^ny^{\prime}_2(1,\mu_n(u),u),\quad \forall n\in\mathbb{N}.\]
Flashka and McLaughlin \cite{flm1976} observed that the  quantities $\{\mu_n,f_n\}_{n\in\mathbb{N}}$ form canonical coordinates of $H^0$, i.e.
\[\{\mu_n,\mu_m\}=\{f_n,f_m\}=0,\quad \{\mu_n,f_m\}=\delta_{n,m},\quad \forall n, m\in \mathbb{N}.\]
Accordingly, the  symplectic form $\omega_2^G$ (see (\ref{omega1})) equals $d\omega_1$, where $\omega_1$ is the 1-form
$\sum_{n\in\mathbb{N}}f_nd\mu_n$. 
Now the KdV action variables are given by the Arnold formula (\ref{arnold1}),
 where $C_n$ is a circle on the invariant torus $\mbox{Iso}(u)$,  corresponding to $\mu_n(u)$. It is  shown in \cite{flm1976}  that 
\[I_n=\frac{2}{\pi}\int_{\lambda_{2n-1}}^{\lambda_{2n}}\lambda\frac{\dot{\triangle}(\lambda)}{\sqrt{\triangle^2(\lambda)-4}}d\lambda,\quad \forall n\in\mathbb{N}.\]
The analytic properties of the functions $u\mapsto I_n$ and of the mapping
$u\mapsto I=(I_1,I_2,\dots)$ were studied later by Kappeler and Korotyaev (see references in
\cite{kjp2003, kor2006} and below). In particular, it was shown that 
 $I_n(u)$, $n\in\mathbb{N}$, are real analytic functions on  $H^0$ of the form 
$I_n=g_n^2+$  {\it higher order terms}, and  $I_n=0$ if and only if  $g_n=0$, 
 see in  \cite{kjp2003}. 
 For any vector $I=(I_1,I_2,\dots)$ with non-negative components we will  denote 
 \begin{equation}\label{tor}
 T_I=\{u(x)\in H^0:  I_n(u)=I_n\quad \forall\, n\}. 
 \end{equation}

  The angle-variables $\varphi^n$ 
   on the finite-gap manifolds $\mathcal{J}^n$  were found in 1970's by Soviet  mathematicians, who constructed them from the Dirichlet spectrum $\{\mu_j(u)\}$ 
   by means of the Abel transform, associated with the Riemann surface of the function $\sqrt{\triangle^2-4}$, see \cite{DMN1976,mar1977,zmn1984}, and see \cite{IM, Dub, Kr, BBE} for the celebrated 
   explicit formulas for angle-variables $\varphi^n$  and for finite-gap
    solutions of KdV in terms of the theta-functions. 
  
  In \cite{kuk1989}  and \cite{kuk2000}, Section~7, 
  the action-angle variables $(I^n, \varphi^n)$ on a finite-gap 
  manifold $\mathcal{J}^n$ and the explicit formulas for solutions of KdV on manifolds 
  $\mathcal{J}^N$, $N\ge n$, from the works \cite{Dub, Kr, BBE}  were used to obtain 
  an analytic symplectic coordinate system $(I^n, \varphi^n, y)$ in the vicinity of
  $\mathcal{J}^n$ in $H^p$. The variable $y$ belongs to a ball in a subspace  $Y\subset H^p$ of 
  co-dimention $2n$, and in the new coordinates  the KdV Hamiltonian (\ref{kdvh})
   reads 
   \begin{equation}\label{SKNF}
 {\cal H}= {\rm const} + h^n(I^n)+\langle A(I^n)y,y\rangle + O(y^3). 
 \end{equation}
 The selfadjoint operator $A(I^n)$ is diagonal in some fixed symplectic basis of $Y$. The nonlinearity $O(y^3)$ defines
 a hamiltonian operator of order one. That is, the KdV's linear operator, which is an operator of order three, 
 mostly transforms to the linear part of the new hamiltonian operator and ``does not spread much" to its nonlinear part.
 This is the crucial property of (\ref{SKNF}).  The
  normal form (\ref{SKNF}) 
   is  instrumental for the purposes of the KAM-theory, see below Section~\ref{skam}.

   McKean and Trubowitz in \cite{mct1976,mct1978} extended the construction  of angles on 
   finite-gap manifolds 
   to the set of all potentials, thus obtaining  angle variables  $\varphi=(\varphi_1,\varphi_2,\dots)$ on the whole space $H^p$, $p\geqslant0$. The angles $(\varphi_k(u),k\ge1)$  are well defined  Gato-analytic functions of $u$ 
   outside the 
   locus
   \begin{equation}\label{game}
   \Game = \{u(x): g_j(u)=0\; {\rm for\  some} \; j\} ,
   \end{equation}
   which is dense in each space $H^p$. The action-map $u\mapsto I$ was not considered 
   in \cite{mct1976,mct1978}, but it may  be shown that  outside $\Game$,  in a certain  weak sense,  the variables 
   $(I,\varphi)$ are KdV's action-angles (see the  next section for a 
   stronger statement).  This  result  is nice and elegant, but it is insufficient 
    to study perturbations of KdV since the transformation to the variables $(I,\varphi)$ is singular at the dense locus 
   $\Game$.

\subsection{Birkhoff coordinates and nonlinear Fourier transform} 
 In a number of publications (see in \cite{kjp2003}),  Kappeler with collaborators proved that the  Birkhoff coordinates $v=
 \{v_n,\;n=\pm1,\pm2,\dots\}$, associated with the action-angles variables $(I,\varphi)$, 
\begin{equation}
v_n=\sqrt{2I_n}\cos(\varphi_n),\quad v_{-n}=\sqrt{2I_n}\sin(\varphi_n),\quad \forall n\in\mathbb{N},
\label{v-variable}
\end{equation}
are analytic on the whole of $H^0$ and define there  a global coordinate system, in which the KdV Hamiltonian (\ref{kdvh})  is a function of  the actions only. This remarkable  result significantly 
specifies 
the normal form (\ref{SKNF}). To 
state it 
exactly, for any $p\in\mathbb{R}$, we introduce the  Hilbert space $h^p$,
\[
\fl
h^p:=\Big\{ v=(\mathbf{v}_1,\mathbf{v}_2,\cdots): |v|_p^2=\sum_{j=1}^{+\infty}(2\pi j)^{2p+1}|\mathbf{v}_j|^2<\infty,\;\mathbf{v}_j=(v_j,v_{-j})^t\in \mathbb{R}^2,\; j\in\mathbb{N}\Big\},\]
and the weighted $l^1$-space $h^p_I$,
\[h^p_I:=\Big\{I=(I_1,\dots)\in\mathbb{R}^{\infty}:|I|^{\sim}_p=2\sum_{j=1}^{+\infty}(2\pi j)^{2p+1}|I_j|<+\infty\Big\}.
\]
Define the mappings
\[
\pi_I:\; h^p\to h_I^p,\quad v\mapsto I=(I_1,I_2,\dots), \;\; {\rm where}\;\;
I_k=\frac{1}{2}| {\mathbf{v}}_k|^2
\quad \forall\, k, \]
\begin{equation*}
\fl
\qquad\eqalign{
{}\qquad\qquad
 \pi_\varphi :\; h^p \to {\mathbb T}^\infty, \quad v\mapsto \varphi&=(\varphi_1,\varphi_2,\dots),\;\;
  {\rm where}\;\; \varphi_k=\arctan(\frac{v_{-k}}{v_k})\;\;
\cr
&  {\rm  if } \;\;  \mathbf{v}_k\ne0,\;\;  {\rm and }\;\; \varphi_k=0\;\; {\rm  if } \;\;  \mathbf{v}_k=0.
}
\end{equation*}
Since $|\pi_I(v)|_p^{\sim}=|v|_p^2$, then 
$\pi_I$ is continuous.  Its image  $h^p_{I+}=\pi_I(h^p)$ is the  positive octant in $h_I^p$. When there is no ambiguity, we  write $I(v)=\pi_I(v)$.

Consider the mapping
 \[\Psi: u(x)\mapsto v=(\mathbf{v}_1,\mathbf{v}_2,\dots),\quad \mathbf{v}_n=(v_n,v_{-n})^t\in\mathbb{R}^2,
 \]  
  where $v_{\pm n}$ are defined by (\ref{v-variable}) and $\{I_n(u)\}$, $\{\varphi_n(u)\}$ are the actions and  angles as in Section~3.2.  Clearly $\pi_I \circ\Psi(u)=I(u)$ and    $\pi_\varphi \circ\Psi(u)=\varphi(u)$.
   Below we refer to $\Psi$ as to   the {\it nonlinear Fourier transform}.

\begin{theorem} (see \cite{kjp2003,kmt2005}) \label{t_kapp}
The mapping $\Psi$ defines an analytical symplectomorphism $\Psi:(H^0,\omega^G_2)\to(h^0,\sum_{k=1}^{\infty}dv_k\wedge dv_{-k})$ with the following properties:
 \begin{enumerate}
 \item[](i) For any $p\in[-1,+\infty)$, it  defines an analytic diffeomorphism $\Psi : H^p \mapsto h^p$.
 \item[](ii) (Percival's identity) If $v=\Psi(u)$, then $|v|_0=||u||_0$.
  \item[](iii) (Normalisation)  The differential  $d\Psi(0)$ is the operator $\sum u_se_s\mapsto v$, where 
   $v_s=|2\pi s|^{-1/2}u_s$ for each $s$. 
 \item[](iv)  The function $\hat{H}(v)=\mathcal{H}(\Psi^{-1}(v))$ has the form $\hat{H}(v)=H_K(I(v))$, where the function $H_K(I)$ is analytic in a suitable neighborhood of the octant $h^1_{I+}$ in $h^1_I$, such that a curve $u\in C^1(0,T;H^0)$ is a solution of  KdV 
   if and only if $v(t)=\Psi(u(t))$ satisfies the equations
  \begin{equation} \dot{\mathbf{v}}_j=\mathbb{J}\frac{\partial H_K}{\partial I_j}(I)\mathbf{v}_j  , \quad\mathbf{v}_j=(v_j,v_{-j})^t\in\mathbb{R}^2,\;j\in \mathbb{N}.
\label{bnf-e}
 \end{equation}
 \end{enumerate}
 \label{Thbnf1}
\end{theorem}
The assertion (iii) normalizes $\Psi$ in the following sense. For any $\theta=(\theta_1,\theta_2,\dots)\in\mathbb{T}^{\infty}$ denote by $\Phi_{\theta}$ the operator
\begin{equation}
\Phi_{\theta}v=v^{\prime}, \quad \mathbf{v}^{\prime}_j=\bar{\Phi}_{\theta_j}\mathbf{v}_j,\quad \forall j\in\mathbb{N},
\label{phi-rotation}
\end{equation}
where $\bar{\Phi}_{\alpha}$ is the rotation of the plane $\mathbb{R}^2$ by the angle $\alpha$. Then $\Phi_{\theta}\circ\Psi$ satisfies all assertions of the theorem except (iii). But the properties (i)-(iv) jointly determine $\Psi$ in a unique way.

The theorem above  can be viewed as a global infinite dimensional version of the Vey
 Theorem~\ref{th-vey}  for KdV, and eq.~(\ref{bnf-e}) -- as a global Birkhoff normal form for KdV. Note that in finite dimension 
 a global Birkhoff normal form exists only for very exceptional integrable  equations, which were found
 during the boom of activity in integrable systems, provoked by the discovery of the method of Lax pair.

\begin{remark}The map $\Psi$ simultaneously transforms  all Hamiltonians of the KdV hierarchy to the Birkhoff normal form. The {\it KdV hierarchy} is a  collection of hamiltonian  functions $\mathcal{J}_l$, $l\geqslant0$, commuting with the KdV Hamiltonian,  and having the  form
\[\mathcal{J}_l(u)= \int
\Big(\frac12 (u^{(l)})^2+ J_{l-1}(u)\Big)dx.\]
Here $J_{-1}=0$ 
 and $J_{l-1}(u),\ l\ge1$, is a polynomial of $u,\dots,u^{(l-1)}$. 
 The functions from the KdV hierarchy form another complete set of KdV integrals.  E.g. see \cite{DMN1976,kjp2003, lax1996}.
\label{remark-laws}
\end{remark}

Properties of the nonlinear Fourier transform $\Psi$ may be specified in two important respects. 
 One of  this specifications  -- the quasilinearity of $\Psi$ -- is   presented in the theorem below. Another one -- its behaviour at infinity -- 
is discussed in the next section.

The  nonlinear Fourier transform $\Psi$ is quasi-linear.   Precisely,
\begin{theorem}
If $m\geqslant0$, then the map $\Psi-d\Psi(0):$ $H^m\to h^{m+1}$ is analytic. 
\label{quasi-l}
\end{theorem}

That is, the non-linear part of $\Psi$ is 1-smoother than its linearisation at the origin. 
 See \cite{kpe2010}  for a local version of this theorem, applicable as well to other integrable infinite-dimensional systems,  and see \cite{KST1, KST2} for the  global result.  We note that the transformation to the normal form (\ref{SKNF}) also 
 is known to be 
 quasi-linear, see  \cite{kuk1989, kuk2000}.

 \begin{problem} Does the mapping $\Psi-d\Psi(0)$ analytically maps $H^m$ to $h^{m+1+\gamma}$ with $\gamma>0$?
 \label{p.smoothing}
 \end{problem}
It is  proved in \cite{kpe2010} that   for $\gamma>1$ the answer to this problem is negative and is
conjectured there that it  also is negative for $\gamma>0$.

\subsection{Behaviour  of $\Psi$  near infinity and large solutions of KdV}

 By the assertion (ii) of Theorem \ref{Thbnf1}, $|\Psi(u)|_0=||u||_0$. It was established by Korotyayev in \cite{kor2006} that higher order norms of $u$ and $v=\Psi(u)$ are related by both-sides polynomial estimates:
 \begin{theorem} For any $m\in\mathbb{N}$, there are polynomials $\mathcal{P}_m(y)$ and $\mathcal{Q}_m(y)$ such that if $u\in H^m$ and $v=\Psi(u)$, then 
 \[|v|_m\leqslant \mathcal{P}_m(||u||_m),\quad||u||_m\leqslant \mathcal{Q}_m(|v|_m).\]
 \label{t.kor}
 \end{theorem}
 The polynomials $\mathcal{P}_m$ and $\mathcal{Q}_m$ are constructed in \cite{kor2006} inductively. From a personal communication of   Korotyayev we know that one can take 
 \begin{equation}
 \mathcal{P}_m(y)=C_my(1+y)^{\frac{2(m+2)}{3}}.
 \label{kor1}
 \end{equation}
 Estimating a potential  $u(x)$ via its actions\footnote{Note that  any $|v|_m^2$ is a linear combinations of the actions $I_j$, $j\geqslant1$.} is more complicated. Corresponding polynomials $\mathcal{Q}_m$ may be chosen of the form
 \begin{equation}
 \mathcal{Q}_m(y)=C^{\prime}_my(1+y)^{a_m},
 \label{kor2}
 \end{equation}
 where $a_1=\frac{5}{2}$, $a_2=3$ and $a_m$ has a factorial growth as $m\to\infty$. 
 \begin{problem} Prove that there exist polynomials $\mathcal{P}^1_m$ and $\mathcal{Q}^1_m$, $m\in\mathbb{N}$, such that for any $u\in H^m$, we have
 \[||d\Psi(u)||_{m,m}\leqslant \mathcal{P}^1_m(||u||_m),\quad ||d\Psi^{-1}(v)||_{m,m}\leqslant \mathcal{Q}^1_m(||v||_m),\]
 where $v=\Psi(u)$. Prove similar polynomial bounds for the norms of higher differentials of $\Psi$ and $\Psi^{-1}$. 
 \label{p.kor}
 \end{problem}
 It seems that to solve the problem a new proof of Theorem \ref{Thbnf1} has to be found (note that  the existing proof  is rather bulky and occupies  half of the book \cite{kjp2003}).
 
 The difficulty in resolving the problem above streams from the fact that $\Psi(u)$ is constructed in terms of spectral characteristics of the Strum-Liouville operator $L_u$, and their dependence on large potentials $u(x)$ is poorly understood. Accordingly, the following question seems to be very complicated:
 \begin{problem}
 Let $u(t,x)=u(t,x;\lambda)$ be a solution of (\ref{kdv1}) such that $u(0,x)=\lambda u_0(x)$, where $u_0\not\equiv0$ is a given smooth  function with zero mean-value, and $\lambda>1$ is a large parameter. The task is to study behaviour of $u(t,x;\lambda)$ when $\lambda\to\infty$.
 \label{p.ll}
 \end{problem}
 Let $u(t,x)$ be as above. Then  $||u(0,x)||_m=\lambda C(m,u_0)$.
 By Theorem~\ref{t.kor}, $|v(0)|_m\leqslant \mathcal{P}_m(\lambda C(m,u_0))$. Since $|v(t)|_m$ is an integral of motion, then using again the theorem we get that 
 \begin{equation*}||u(t)||_m\leqslant \mathcal{Q}_m\Big(\mathcal{P}_m(\lambda C(m,u_0))\Big).
 \end{equation*}
 In particular, by (\ref{kor1}) and (\ref{kor2}) we have
 $\ 
 ||u(t)||_1\leqslant C(1+\lambda)^{21/2}.
 $
 A lower bound for the Sobolev norms comes from the fact that $||u(t)||_0$ is an integral of motion. So 
 \[||u(t)||_m\geqslant ||u(0)||_0=\lambda C(0,u_0).\]
  We have demonstrated:
 \begin{proposition}
 Let $u(t,x)$ be a solution of (\ref{kdv1}) such that 
 $u(0,x)=\lambda u_0(x)$,  where $ 0\not\equiv
   u_0\in C^{\infty}\cap H^0$  and $\lambda\geqslant1$.
 Then for $m\ge1$ we have 
 \begin{equation}
 \fl
 \hspace{40pt}
 1+(\lambda C(m,u_0))^{A_m}\geqslant\limsup_{t\to\infty}||u(t)||_m\geqslant\liminf_{t\to\infty}||u(t)||_m\geqslant c(u_0) \lambda,\quad
 \label{star-3}
 \end{equation}
 \end{proposition}
 for a suitable $A_m>1$.
 E.g. $A_1=21/2$.
 \smallskip
 
 The third estimate in (\ref{star-3})  is  optimal up to a constant factor as 
 \[\liminf_{t\to\infty}||u(t)||_m\leqslant \lambda C(m,u_0),\]
 since the curve $u(t)$ is almost periodic. The first estimate with the exponent $A_m$ which follows  from
 Theorem~\ref{t.kor} 
  certainly is not optimal. But the assertion   that $\limsup_{t\to\infty}||u(t)||_m$ grows with $\lambda$ as $\lambda^{A_m}$, where $A_m$ goes to infinity with $m$, 
   is correct. It follows from our next result:
 \begin{theorem}
 Let $k\geqslant4$. Then there exists $\alpha>0$ and, for any $\lambda>1$ there exists $t_{*}=t_{*}(u_0,\lambda)$ such that 
 \begin{equation}
 ||u(t_{*})||_k\geqslant c_{u_0}^{\prime}\lambda^{1+\alpha k}.
 \label{gafa}
 \end{equation}
 \label{t.gafa}
 \end{theorem}
 In \cite{kgafa99} (see there Theorem 3 and Appendix 2), the theorem is proved for a class of non-linear Schr\"odingier equations which includes the defocusing Zakharov-Shabat equation. The proof applies to KdV.
 See \cite{bir2004}, where the argument is applied to the multidimensional 
 Burgers equation,  similar to KdV for the proof of this result. 
 We  mention that (in difference with majority of 
  results in this work) the assertion of  Theorem \ref{t.gafa} remains true for other boundary conditions.
 
 Problem \ref{p.ll} may be scaled as the non-dispersive limit for KdV. Indeed, let us substitute $u=\lambda w$ and pass to fast time $\tau=\lambda t$. Then the function $w(\tau,x;\lambda)$ satisfies 
 \begin{equation}
 w_{\tau}+\lambda^{-1}w_{xxx}-6ww_x=0,\quad w(0,x)=u_0(x),
 \label{LL-1}
 \end{equation}
 and we are interested in $w(\lambda t,x;\lambda)$ when $\lambda\to\infty$. For $\lambda=\infty$ the equation above becomes the Hopf equation. Since $u_0(x)$ is a periodic non-constant function, then the solution of (\ref{LL-1})$_{\lambda=\infty}$ developes a shock at time $\tau_{*}$, $0<\tau_{*}<\infty$. Accordingly, the elementary perturbation theory allows to study solutions of (\ref{LL-1}) when $\lambda\to\infty$ for $\tau<\tau_{*}$, but not for $\tau\geqslant\tau_{*}$. The problem to study this limit for $\tau\geqslant \tau_{*}$ is addressed by the Lax-Levermore theory (mostly for the case when $x\in\mathbb{R}$ and $u_0(x)$ vanishes at infinity). There is vast literature on this subject, e.g. see \cite{LL,dvz1997} and references  in \cite{dvz1997}. The existing Lax-Levermore theory does not allow to study solutions $w(\tau,x)$ for $\tau\sim\lambda^{-1}$,  as is required by Problem \ref{p.ll}.

 \subsection{Properties of frequency map}
  Let us denote 
 \begin{equation}W(I)=(W_1(I),\;W_2(I),\dots),\quad W_i(I)=\frac{\partial H_K}{\partial I_i},\quad i\in\mathbb{N}.
 \label{frequency}
 \end{equation}
 This is the {\it frequency map for  KdV.}  By Theorem \ref{Thbnf1} each its component is an analytic function, defined in the vicinity of $h^1_{I+}$ in $h^1_I$.
\begin{lem} 
a) For $i,j\ge1$ we have 
$\partial^2W(0) /\partial I_i \partial I_j=-6\delta_{i,j}$.

b) For any $n\in\mathbb{N}$, if $I_{n+1}=I_{n+2}=\dots=0$, then
\[\det\Big(\big(\frac{\partial W_i}{\partial I_j}\big)_{1\leqslant i,j\leqslant n}\Big)\neq0.\]
\label{lem-nd}
\end{lem}

For a) see  \cite{bok1991, kjp2003,kuk2000}. For  a proof of b) and references to the original works
of Krichever and Bikbaev-Kuksin  see  Section~3.3 of \cite{kuk2000}. 

Let $l^{\infty}_i$, $i\in\mathbb{Z}$, be the Banach spaces of all real sequences $l=(l_1,l_2,\dots)$ with norms
\[|l|^{\infty}_i=\sup_{n\geqslant 1} n^i|l_n|<\infty.\]
Denote $\boldsymbol{\kappa}=(\kappa_n)_{n\in\mathbb{N}},$ where $\kappa_n=(2\pi n)^3$. For the following result see \cite{kjp2003}, Theorem~15.4.

\begin{lem} The normalized frequency map
$\  I\mapsto  W(I)-{\boldsymbol{\kappa}}$
is real analytic as a mapping from $h^1$ to $l^{\infty}_{-1}$.
\label{lem-reg}
\end{lem}

 From these two lemmata we known that the Hamiltonian $H_K(I)$ of KdV is non-degenerated in the sense of Kolmogorov and its nonlinear part is more regular than its linear part.  These properties are very important  to study  perturbations of KdV.

\subsection{Convexity of  Hamiltonian $H_K(I)$}
By Theorem \ref{Thbnf1},  the dynamics of KdV is determined by  the Hamiltonian $H_K(I)$.  To understand the  properties of the latter is an important step toward  the study of perturbations of KdV.

Denote by $P_j$ the moments of the actions, given by 
\[P_j=\sum_{i\geqslant1}(2\pi n)^jI_i,\quad j\in\mathbb{Z}.\]
Due to Theorem \ref{Thbnf1}, the linear part of $H_K(I)$ at the origin $dH_K(0)(I)$  equals to $\frac{1}{2}\int_{\mathbb{T}}(\frac{\partial}{\partial x}(\Psi^{-1}v))^2dx=P_3$. So we can write $H_K(I)$ as 
\[H_K(I)=P_3(I)-V(I), \quad V(I)=\mathcal{O}(|I|^2_1).\]
(The minus-sign here is convenient since, as we will see,  $V(I)$ is non-negative.) For any $N\geqslant 1$, denote $\tilde{l}^N\subset l^2$ the $N$-dimensional subspace
\[\tilde{l}^N=\{l=(l_1,\dots)\in\mathbb{R}^{\infty}: l_n=0,\;\forall n>N\},\]
and set $\tilde{l}^{\infty}=\cup_{N\in\mathbb{N}}\tilde{l}^N$. Clearly the function $V$ is analytic on each octant $\tilde{l}^N_+$. So it is Gato-analytic on the octant 
 $\tilde{l}^{\infty}_+$. That is, it is analytic on every interval 
 $\{(ta+(1-t)c)\in\tilde{l}_+^{\infty}: 0\le t\le1 \}$, where $a,c\in\tilde{l}^{\infty}_+$.

 By Lemma~\ref{lem-nd}~a), 
 $d^2V(0)(I)=6||I||_2^2$. This suggests that the Hilbert space $l^2$ rather than the Banach space $h_I^1$ (which is contained in $l^2$) is a distinguished phase space for the Hamiltonian $H_K(I)$. This guess is justified by the following result:
\begin{theorem} (see  \cite{kok2011}). 
(i) The function $V: \tilde{l}_+^{\infty}\to \mathbb{R}$ extends to a non-negative continuous function on the $l^2$-octant $l^2_+$, such that $V(I)=0$ for some $I\in l^2_+$ iff $I=0$. Moreover 
$\ 0\leqslant V(I)\leqslant 8P_1P_{-1}.$

(ii) For any $I\in l^2_+$, the following estimates hold true:
$$
\frac{\pi}{10}\frac{||I||_2^2}{1+2P_{-1}^{1/2}}\leqslant V(I)\leqslant (8^3(1+P_{-1}^{1/2})^{1/2}P_{-1}^2+6\pi  e^{P_{-1}^{1/2}/2}||I||_2)||I||_2.
$$
(iii) The function $V(I)$ is convex on $l_+^2$.
\label{th-convex}
\end{theorem}

Note that the assertion (iii) follows from (i) and Lemma \ref{lem-nd}. Indeed, since $V(I)$  is analytic on $\tilde{l}^N_+$, then Lemma \ref{lem-nd}  assures that the Hessian $\{\frac{\partial^2V}{\partial I_i\partial I_j}\}_{1\leqslant i,j\leqslant N}$ is positive definite on $\tilde{ l}_+^N$. Thus $V$ is convex on $\tilde{l}_+^N$, for each $N\in\mathbb{N}$.  Then the assertion (iii) is deduced from the fact the $\tilde{l}_+^{\infty}=\cup_{N\in\mathbb{N}}\tilde{l}_+^N$ is dense in $l^2_+$, where $V(I)$ is continuous.

\smallskip
\smallskip

Assertion ii) of the theorem shows  that $l^2$ is the biggest Banach space on which $V(I)$ is continuous\footnote{That is, if $V(I)$ continuously extends to a Banach space $\mathcal{B}$ of the sequences $(I_1,I_2,\dots)$, then $\mathcal{B}$ may be continuously embedded in $l^2$.}. Jointly with the convexity of $V(I)$ on $l^2$, this hints that $l^2$ is the natural space to study the long-time dynamics of actions $I(u(t))$ for solutions $u(t)$ of a perturbed KdV.

This theorem  and the analyticity  of the KdV Hamiltonian (\ref{kdvh}) do not leave much doubts that $V(I)$ is analytic on $l^2_+$. If so, then by Lemma~\ref{lem-nd}, this function is strictly convex in a neighbourhood of the origin in $l^2$. Most likely, it is strictly convex everywhere on $l^2_+$.

 In difference with $V(I)$, the total Hamiltonian $H_K(I)$ is not continuous on $l^2_+$ since its linear part $P_3(I)$ is there an unbounded linear functional. But  $P_3(I)$ contributes to equation (\ref{bnf-e}) the linear rotation
$\dot{\mathbf{v}}_k=\mathbb{J}(2\pi n)^3 \mathbf{v}_k$, $ k\in\mathbb{N}$.
So the properties of equation (\ref{bnf-e}) essentially are determined by the component $V(I)$ of the Hamiltonan.  Note that since $P_3(I)$ is a bounded linear functional on the space $h^1\subset l^2$, then the Hamiltonian $H_K(I)$ is concave in $h^1$.
\begin{problem}
Is it true that the function $V(I)-3|I|_{l^2}^2$ extends analytically (or continuously) to a space, bigger than $l^2$?
\end{problem}


\section{ Perturbations of KdV}
In the theory of integrable systems in  finite dimension, there are two types of  perturbative results concerning  long-time   stability of  solutions.  The first type is  the KAM theory. Roughly, it says that among a family of invariant tori of the unperturbed system,
given by the  Liouville-Arnold-Jost theorem,  there exists (under generic assumptions) a large set of tori which  survive under  sufficiently small hamiltonian  perturbations, deforming only slightly. In particular, the perturbed system admits plenty of quasi-periodic solutions 
\cite{AKN, arn1963}.  Results of the second type are obtained by the techniques of averaging which applies to a larger class of dynamical systems, characterized by the existence of fast and slow variables. This method has a much longer history which dates back to the epoch  of Lagrange and Laplace, who  applied it to the  problems of celestial mechanics, without proper justifications. Only in the last fifty years  rigorous mathematical justification  of the principle has been  obtained, see
 in \cite{ nei1975, AKN, lom1988}. If the unperturbed system is hamiltonian  integrable, then for  the slow-fast variables one can  choose the action-angle variables. The averaging theorems say that under appropriate assumptions, the action variables, calculated for solutions  of the perturbed system, can be well approximated by  solutions of a suitable averaged vector field, over an extended time interval. If the perturbation is hamiltonian, then the averaged 
 vector filed vanishes.  The strongest result in this direction is due to Nekhoroshev \cite{nek1972,loc1992}, who proved that in the  hamiltonian  case  the action variables vary just a bit over  exponentially long time intervals.

 Concerning  instability of   solutions,  also two types of phenomena are known. 
  One is called the Arnold diffusion. In \cite{arn1964}, Arnold observed that despite  most of  the phase space of near integrable hamiltonian   systems with more than two degrees of freedom  is foliated by invariant  KAM tori, still there can exist solutions such that their actions admit increments of order one  during  sufficiently long time. Arnold conjectured that this phenomenon is generic. Though there are many developments in this direction in the last ten years, the mechanism  of the  Arnold diffusion  still is  far from being  well understood. Another instability mechanism  is known as the  capture in resonance.   The essence of this phenomenon is that a solution of a perturbed system 
   reaches a resonant zone and begins drifting along it in such a way  that the resonance condition approximately holds. Therefore, solutions of the original perturbed  system and the averaged one diverge by a quantity of order one on  a time interval of order $\epsilon^{-1}$ (see e.g. \cite{nei2005}).  
   
   Now return to PDEs. Two types of perturbations of the KdV equation (\ref{kdv1}) have been 
   considered: when the boundary condition is perturbed but the equation is not, and 
   other way round. Problems of the first type lie outside the scope of our work, and very
   few results (if any) are proved there rigorously, see \cite{Kr}  for discussion and some 
   related statements. Problems of the second type are much closer to the finite-dimensional
   situation and  they are discussed below. 
Several  attempts were made to establish stability results for perturbed integrable PDEs, e.g. for perturbed KdV, analogous to those in finite dimension.   Among them, the KAM theory was  the most successful \cite{Craig}. The first results in this direction  are  due to Kuksin \cite{kuk1987, kuk1989} and Wayne \cite{way1990}. Despite there are no  rigorously proven  instability results for the perturbed KdV,  we mention our believe that to study the instability, Theorem \ref{Thbnf1}, \ref{quasi-l} and \ref{th-convex} should be important. 

\subsection{KAM theorem for perturbed KdV}\label{skam}
Consider  the  hamiltonian  perturbation of KdV, corresponding to a Hamiltonian $H_{\epsilon}=\mathcal{H}(u)+\epsilon F(u)$:
\begin{equation}
\dot{u}+u_{xxx}-6uu_x-\epsilon \frac{\partial}{\partial x}\nabla F(u)=0,\quad u(0)\in H^1.\label{pkdvh1}
\end{equation}
Here $F(u)$ is an analytic functional on $H^1$ and $\nabla F$ is its $L^2$-gradient in $u$.  Let $n\in\mathbb{N}$, and $\Gamma\subset \mathbb{R}_+^n$ be a compact set of positive Lebesgue measure. Consider a family of the $n$-gap tori:
\[\mathcal{T}_{\Gamma}=\cup_{I\in\Gamma}\mathcal{T}^n_I\subset\mathcal{J}^n,\quad
 \mathcal{T}^n_I=T_{(I,0,\dots)},
\]
where $I=(I_1,\dots,I_n)$,  see (\ref{finite-gap}) and (\ref{tor}). It turns out that most of them
 persist as invariant tori of the perturbed equation (\ref{pkdvh1}):
\begin{theorem}
For some 
 $M\geqslant1$, assume that the Hamiltonian $F$  analytically extends to a complex 
   neighbourhood $U^c$ of $\mathcal{T}_{\Gamma}$ in $H^M \oplus\mathbb C$ and satisfies there the regularity condition
\[\nabla F: U^c
\to H^M \otimes\mathbb C
,\quad \sup_{u\in U^c} (|F(u)|+ ||\nabla F(u)||_M)\leqslant1.\]
Then, there exists an $\epsilon_0>0$ 
 and for $\epsilon<\epsilon_0$ there exist

(i) a nonempty Cantor set $\Gamma_{\epsilon}\subset \Gamma$ with $\mbox{mes}(\Gamma\setminus\Gamma_{\epsilon})\to 0$ as $\epsilon \to 0$;

(ii) a Lipshcitiz mapping 
$\ \Xi: \mathbb{T}^n\times\Gamma_{\epsilon}\to U\cap H^M$, 
such that its restriction to each torus $\mathbb{T}^n\times I$, $I\in\Gamma_{\epsilon}$, is an analytical embedding;

(iii) a Lipschitz map $\chi: \Gamma_{\epsilon}\to \mathbb{R}^n$, $|\chi-\mbox{Id}|\leqslant Const\cdot \epsilon$,  

\noindent 
such that for every $(\varphi,I)\in \mathbb{T}^n\times\Gamma_{\epsilon}$, the curve $u(t)=\Xi(\varphi+\chi(I)t,I)$ is an quasi-periodic solution of (\ref{pkdvh1}) winding around the invariant torus $\Xi(\mathbb{T}^n\times\{I\})$. Moreover, these solutions are  linearly stable.
\label{kdv-kam}
\end{theorem}
\noindent 
{\it Proof:}\quad In the coordinates $v$ as in  Theorem \ref{Thbnf1}, the Hamiltonian $H_{\epsilon}$ becomes $H_K(I)+\epsilon F(v)$. For any $I_0\in h_I^M$, using Taylor's formula,  we write 
\begin{equation}
\fl
\qquad\eqalign{H_K(I_0+I)&=H_K(I_0)+\sum_{i\geqslant1}\frac{\partial H_K}{\partial I_i}(I_0)I_i+\int_0^1(1-t)\sum\frac{\partial^2 H_K}{\partial I_i\partial I_j}I_iI_jdt\cr
&:=const+\sum_{i\geqslant1}W_i(I_0)I_i+Q(I_0,I),}
\label{h-exp}
\end{equation}
where $W$ is the frequency map, see (\ref{frequency}).
Now we introduce the symplectic polar coordinates around the tori in the family $\mathcal{T}_{\Gamma}$.
Namely,  for each $\xi_0\in\Gamma$, we set
\[\cases{v_i=\sqrt{\xi_0+y_i}\cos\varphi,\quad v_{-i}=\sqrt{\xi_0+y_i}\sin\varphi, &$1\leqslant i\leqslant n,$\\ b_i=v_i, \quad b_{-i}=v_{-i}, & $i\geqslant n+1.$\\}
\]
Denote $b=(b_{n+1},b_{-n-1},\dots)$. The transformation above  is real analytic and symplectic on 
\[D(s,r)=\{|{\rm Im}\,\varphi |<s\}\times\{|y|<r^2\}\times\{|b|_M<r\},
\]
for all $s>0$ and $r>0$ small enough. Using the expansion (\ref{h-exp}), setting $I_0=(\xi_0,0)$ and neglecting an irrelevant constant we see that the integrable Hamiltonian in new coordinates is given by 
\[H_K=N+Q=N(y,\xi_0,b)+Q(y,\xi_0,b),\]
where $N=\sum_{1\leqslant i\leqslant n}W_i(\xi_0)y_i+\frac{1}{2}\sum_{i\geqslant 1} W_{n+i}(b_{n+i}^2+b_{-n-i}^2),$
with $I_i=y_i$ for $1\leqslant i\leqslant n$, and $2I_i=b_i^2+b_{-i}^2$ for $i>n$. Then the whole Hamiltonian of the perturbed equation (\ref{pkdvh1}) can be written as
\begin{equation}H_{\epsilon}=N+Q+\epsilon F.
\label{pkdvhnf}
\end{equation}
We consider the new perturbation term $P=Q+\epsilon F$. By Lemma \ref{lem-reg} and the regularity assumption, if $r^2=\sqrt{\epsilon}$, then 
\[\sup_{D(s,r)}||X_P||_{M-1}\leqslant c\sqrt{\epsilon}.
\]
Using the non-degeneracy  Lemma \ref{lem-nd}, we can take the vector 
$\{W_i(\xi_0), 1\le i\le n\}$ for a free  $n$-dimensional parameter of the problem and  apply an abstract KAM theorem (see Theorem 8.3 in \cite{kuk2000} and  Theorem 18.1 in \cite{kjp2003}) to obtain the statements of the theorem.
For a complete proof, see \cite{kuk2000,kjp2003}.
\quad $\square$
\smallskip

\noindent
{\it Example.}
 The theorem applies if (\ref{pkdvh1}) is the hamiltonian  PDE with the local Hamiltonian
$$
\mathcal{H}(u)+ \epsilon \int_{\mathbb{T}} g(u(x),x)\,dx
=\int_{\mathbb{T}}\big(\frac{u_x^2}{2}+u^3+  \epsilon g(u(x) ,x) \big)\,
dx,
$$
where $g(u,x)$ is a smooth function, periodic in $x$ and analytic in $u$. 
In this case in  (\ref{pkdvh1})  we have 
$\frac{\partial}{\partial x}
\nabla F(u(\cdot))(x)=\frac{\partial}{\partial x} g_u(u(x),x)$.  $\quad\square$
\medskip

It is easy to see that for  the proof of  Theorem \ref{kdv-kam}, explained above, the normal
 form (\ref{SKNF})   (weaker and much simpler 
  than that in Theorem \ref{Thbnf1})   is sufficient, see \cite{kuk1989,kuk2000}.  
  Normal forms, similar to (\ref{SKNF}),   exist for  integrable PDEs  for which a normal form as in Theorem \ref{Thbnf1} does not hold,  e.g, for the \mbox{Sine-Gordon} equations, see \cite{kuk2000}.
  \medskip

\noindent{\bf Remark}. 
Recently Berti and Baldi announced a KAM-theorem, similar to Theorem~\ref{kdv-kam}, which 
applies to perturbations of KdV,  
given by operators of order higher than one.
\medskip

Theorem \ref{kdv-kam} allows to perturb a set of KdV-solutions which form a null-set for any reasonable measure in the function space (in difference with the finite-dimensional KAM~theory which insures the persistence of a set of almost-periodic solutions  which  occupy the phase space up to a set of measure~$\lesssim \epsilon^{\gamma},$ $\gamma>0$). It gives rise to a natural question:

\begin{problem}
  Do typical tori $T_I$ as in (\ref{tor}), 
   where $I=(I_1,I_2,\dots)\in h^p_I$, $p<\infty$, persist in the perturbed equation (\ref{pkdvh1})? If they do not, what happens to them?
 \end{problem}
 
 Though there are KAM-theorems for perturbations of infinite-dimensional invariant tori, e.g. see \cite{pos1990,bou1996}, they are not applicable  to the problem above since, firstly, those works do not apply to KdV due to the strong nonlinear effects and long-range coupling between the modes and, secondly, they only treat invariant tori corresponding to the actions $I_1,I_2, \dots$, decaying  very fast (faster than exponentially).
 
\subsection{ Averaging  for perturbed KdV}
Compare  to  KAM, averaging type theorems for perturbed  KdV are more recent and  less developed. Their stochastic versions, which we discuss in Section~\ref{s_stoch},
 are significantly stronger than corresponding 
deterministic statements in Section~\ref{s_determ}. 
We will explain reasons for that  a bit later. Let us start with 
 the `easiest'  case, where KdV is stabilized by small dissipation:
\begin{equation}
\dot{u}+u_{xxx}-6uu_x=\epsilon u_{xx}, \quad u(0)\in H^3. 
\label{pkdvp1}
\end{equation}
A simple calculation shows that a solution $u(t)$ satisfies 
\[||u(t)||_0\leqslant e^{-\epsilon t}||u(0)||_0.\]
So $u(t)$ becomes negligible  for $t\gg\epsilon^{-1}$. But what happens during  time intervals of order $\epsilon^{-1}$?  Let us  pass to the slow time $\tau=\epsilon t$ and apply to equation (\ref{pkdvp1}) 
 the nonlinear Fourier transform $\Psi$,  denoting  for $k=1,2,\dots$, 
$\Psi_k(u)=\mathbf{v}_k$ 
if $\Psi(u)=v=(\mathbf{v}_1,\dots)$:
\begin{equation}
\fl
\frac{d}{d\tau}\mathbf{v}_k=\epsilon^{-1}\mathbb{J}W_k(I)\mathbf{v}_k+d\Psi_k(\Psi^{-1}(v))(\triangle \Psi^{-1}(v))=: \epsilon^{-1}\mathbb{J}W_k(I)\mathbf{v}_k+ P_k(v),
\;\; k\in\mathbb{N}.
\label{pkdvp1-v}
\end{equation}
Since $I_k=\frac{1}{2}|\mathbf{v}_k|^2$ is an integral of motion  for  KdV, then
\begin{equation}
\frac{d}{d\tau}I_k=\big(
P_k(v),
\mathbf{v}_k\big):=F_k(v),\quad k\in\mathbb{N},
\label{pkdvp1-i}
\end{equation}
where $(\cdot,\cdot)$ stands for the Euclidian scalar product in $\mathbb{R}^2$.  
 Using (\ref{pkdvp1-v}) we get
\begin{equation}
\frac{d}{d\tau}\varphi_k=\epsilon^{-1}W_k(I)+\langle \mbox{term of order 1}\rangle,\;\mbox{if}\; \mathbf{v}_k\neq0,\; k\in\mathbb{N}.
\label{pkdvp1-a}
\end{equation}
We have written equation (\ref{pkdvp1}) in the action-angle variables  $(I,\varphi)$.
Consider the   {\it averaged equation for  actions:}
\begin{equation}
\fl
\frac{d}{d\tau}J_k=\langle F_k\rangle (J),\quad \langle F_k\rangle(J) =\int_{\mathbb{T}^{\infty}}F_k(J,\varphi)d\varphi, \quad k\in\mathbb{N},\quad J(0)=I(u(0))\,,
\label{pkdvp1-av}
\end{equation}
where $F_k(I,\varphi)=F_k(v(I,\varphi))$, $k\in\mathbb{N}$, and $d\varphi$ is the Haar measure on the infinite dimensional torus $\mathbb{T}^{\infty}$.  The main problem of the averaging theory  is to see if the following  holds true:

\smallskip
\noindent{\bf Averaging principle}: Fix any $T>0$. Let $(I(\tau),\varphi(\tau))$ be a solution of (\ref{pkdvp1-i}), (\ref{pkdvp1-a}),  and $J(\tau)$ be  a solution of (\ref{pkdvp1-av}). Then  (either for all or, for  'typical' initial data $u(0)$)  we have
\[||I(\tau)-J(\tau)||\leqslant \rho(\epsilon),\quad \forall\, 0\leqslant \tau\leqslant T,\]
where $||\cdot||$ is a suitable norm, and $\rho(\epsilon)\to 0$ with $\epsilon\to0$.

\smallskip

The main obstacles to prove this  for the perturbed equation (\ref{pkdvp1}) are the following:
\begin{enumerate}
\item[](1) The KdV-dynamics on some tori is resonant.

\item[](2) The well-posedness of the averaged equation (\ref{pkdvp1-av}) is unknown.

\item[](3)   Equations (\ref{pkdvp1-a}) are singular at the locus $\Game$ (see (\ref{game})). 
\end{enumerate}

To handle the third difficulty  observe that eq.  (\ref{pkdvp1-a}) with a specific $k$ 
is  singular when $I_k$ is small. But then  $v_k$ is small,  the $k$-th mode does not
affect much the dynamics, so the equation for $\varphi_k$ may be excluded from the considereation.
Concerning the second  difficulty, in \cite{kuk2010}  the second author of this work established  that the averaged equation (\ref{pkdvp1-av}) may be lifted  to  a regular system in the space $h^p$, which is well posed,  at least locally. More specifically,  note that equation (\ref{pkdvp1-av}) may be written as follows:
\begin{equation}
\frac{d}{d\tau}J_k=\langle F_k\rangle=\int_{\mathbb{T}^{\infty}}(\bar{\Phi}_{\theta_k}\mathbf{v}_k, \;P_k(\Phi_{\theta}v))d\theta=(\mathbf{v}_k,\;R_k(v)),
\label{effective-r}
\end{equation}
$$
R_k=\int_{\mathbb{T}^{\infty}}\bar\Phi_{-\theta_k}P_k(\Phi_{\theta}v)d\theta,
$$
where the maps $\Phi_{\theta}$ and $\bar{\Phi}_{\theta_k}$ are defined in (\ref{phi-rotation}).  
Now consider equation
\begin{equation}
\frac{d}{d\tau} v=R(v).
\label{effective-e}
\end{equation}
Then relation (\ref{effective-r}) implies:
\begin{lem}
If $v(\tau)$ satisfies (\ref{effective-e}), then $I(v(\tau))$ satisfies (\ref{pkdvp1-av}).
\end{lem}
Equation (\ref{effective-e}) is called the {\it  effective equation} for the perturbed KdV equation (\ref{pkdvp1}). It  is rotation-invariant:  if $v(\tau)$ is a solution of (\ref{effective-e}), then for each $\theta\in\mathbb{T}^{\infty}$, $\Phi_{\theta}(v(\tau))$ also is a solution.
 Since  the  map $\Psi$ is quasilinear by Theorem \ref{quasi-l}, we may 
   write $R(v)$ more explicitly. 
  Namely, denote by $\hat{\triangle}$ the Fourier-image of the Laplacian,  $\hat{\triangle}=\mbox{diag}\{-k^2,\;k\in\mathbb{N}\}$, and set  
 \[L=d\Psi(0),\; \;\Psi_0=\Psi-L,\; \; G=\Psi^{-1}=L^{-1}+G_0.\]
  Then $G_0: h^s\to H^{s+1}$ is analytic for any $s\geqslant0$, and direct calculation shows that 
\[R(v)=\hat{\triangle}v+R_0(v),\]
where \[R_0(v)=\int_{\mathbb{T}^{\infty}}[\Phi_{-\theta}L\triangle(G_0\Phi_{\theta}v)+\Phi_{-\theta}d\Psi_0(G\Phi_{\theta}v)\triangle(G\Phi_{\theta}v)]d\theta.\]
Hence $R_0(v)$  is an operator of order one and  the effective equation (\ref{effective-e}) is a Fourier transform of a quasi-linear heat equation with a non-local nonlinearity of first order. Such equations are    locally well posed. Due to the direct relation between the effective equation and the averaged equation, the former  can be used   to study the latter. 

The first difficulty is serious. Sometime it cannot be overcome, and then the averaging fails. A way 
to handle it is discussed in the next section.

\subsection{ Stochastic averaging. }\label{s_stoch}
A way to 
 handle the first obstacle -- the   resonant tori -- 
  is  to add to the perturbed equation (\ref{pkdvp1}) a random force which would shake  solutions 
 $u(t)$ off a resonant torus (as well as off any other invariant torus $T_I$). 
  So let us consider a randomly perturbed KdV:
  \begin{equation}\label{pkdvp-r}
\dot{u}+u_{xxx}-6uu_x=\epsilon u_{xx}+\sqrt{\epsilon}\, \eta(t,x),\quad u(0)=u_0\in C^{\infty}\cap H^0,
 \end{equation}
 $$
 \eta(t,x)=\frac{\partial}{\partial t}\sum_{j\in\mathbb{Z}_0}b_j\beta_j(t)e_j(x).
$$
Here  $\mathbb{Z}_0$ is the set of all non-zero integers, $\{e_j(x)\}_{j\in\mathbb{Z}_0}$ is the basis  (\ref{basis1}) and 
\begin{itemize}
\item all $b_j>0$ and decay fast when $|j|\to \infty$.
\item $\{\beta_j(t)\}_{j\in\mathbb{Z}_0}$ are independent standard Wiener processes (so $\beta_j(t)=\beta_j^{\omega}(t)$, where $\omega$ is a point in a probability space $(\Omega,\mathcal{F},\mathbf{P})$).
\end{itemize}
The scaling factor $\sqrt{\epsilon}$ in the r.h.s is natural since only with this scaling do solutions of equation (\ref{pkdvp-r}) remain of order one when $t\to\infty$ and $\epsilon\to0$.  To simplify formulas we assume that $b_j=b_{-j}$ for all $j$.

In \cite{kpi2008,kuk2010},  Kuksin and  Piatnitski justified the averaging principle for the stochastic equation (\ref{pkdvp-r}).
To explain their result we  pass to the slow time $\tau=\epsilon t$ and use It\^o's formula (e.g. see in \cite{kar1991}) to  write the corresponding equation for the vector of actions 
$I(u(\epsilon^{-1}\tau))=I^{\omega}(\tau)$:
\begin{equation}
\frac{dI_k}{d\tau}= F_k(I,\varphi)+K_k(I,\varphi)+\sum_j G_k^j(I,\varphi)\frac{\partial}{\partial \tau}\beta_j(\tau),\quad k\ge1. 
\label{pkdvpr-i}
\end{equation}
Here $F$ is defined as in (\ref{pkdvp1-i}), $K$ is the Ito term 
$$K_k=\frac{1}{2}\sum_{j\in\mathbb{Z}_0}b_j^2\big((d^2\Psi_k(u)[e_j,e_j],\mathbf{v}_k)+|d\Psi_k(u)e_j|^2\big),
$$
  and $G$ is the dispersion matrix, $G_k^j=b_jd\Psi_k(u)e_j$.

Let us average the equation above:
\begin{equation}
\frac{dJ}{d\tau}=\langle F\rangle (J)+\langle K\rangle(J)+\sum_j\langle G^j\rangle(J)\frac{\partial}{\partial \tau}\beta_j(\tau).
\label{pkdvpr-av}
\end{equation}
Here  $\langle F\rangle$ are the same as in  (\ref{pkdvp1-av}), $\langle K\rangle$ is the average of $K$ and $\langle G^j\rangle(J)$, $j\in\mathbb{Z}_0$, are  column-vectors, forming an  infinite matrix $\langle G\rangle (J)$. The latter is defined as a square root of the averaged diffusion matrix
$\int_{\mathbb{T}^{\infty}}G(J,\varphi)G^{T}(J,\varphi)d\varphi,$
where $G(J,\varphi)$ is formed by the columns $G^j(J,\varphi)$.  Similar to Section 4.2, equation (\ref{pkdvp-r}) also admits an effective equation of the form
\begin{equation*}
\frac{d\mathbf{v}_k}{d\tau}=-k^2\mathbf{v}_k+R^{\prime}_k(v)+\sum_{j}(R^{\prime\prime})^j_k(v)\frac{\partial}{\partial t}\beta_j(\tau),
\qquad k\ge1, 
\end{equation*}
where $R^{\prime}(v)$ is an operator of first order, $R^{\prime\prime}(v)$ is a Hilbert-Schmidt matrix, which is an analytic function of $v$, and $\{\beta_j(\cdot),  j\geqslant1\}$ are standard independent Wiener processes. This is a quasilinear stochastic 
heat equation with a non-local nonlinearity, written in the Fourier
coordinates. It is well posed in the spaces $h^p, p\ge1$. 

Similar to above, if $v(\tau)$ is a solution of (\ref{effective-r}), then $I(v(\tau))$ is a weak solution of (\ref{pkdvpr-av}). 
See \cite{kpi2008, kuk2010}  for details. 

We recall (e.g. see \cite{kar1991})  that a random process $J=J^{\omega}(\tau)$ is a
{\it  weak solution (in the sense of stochastic analysis)} of equation (\ref{pkdvpr-av}), if for almost every $\omega$ it satisfies the integrated version of equation (\ref{pkdvpr-av}), where the processes $\beta_j$'s are replaced by some other independent standard  Brownian motions $\hat{\beta}_j$'s: 
$$
J^{\omega}(\tau)=\int_0^{\tau}(\langle F\rangle+\langle K\rangle)
(J^{\omega}(s))ds+\int_0^{\tau}\sum_j\langle G^j\rangle (J^{\omega}(s))d\hat{\beta}^{\omega}_j(s),\quad \forall \tau\in[0,T].
$$

Fix any $T>0$. Let $u^{\epsilon}(t)$, $0\leqslant t\leqslant \epsilon^{-1}T$, be a solution of (\ref{pkdvp-r}). Introduce slow time $\tau=\epsilon t$ and denote $I^{\epsilon}(\tau)=I(u^{\epsilon}(\epsilon^{-1}\tau))$. Consider the distribution of this random process. This is a measure in the space $C([0,T],h^p_I)$. We assume $p\geqslant 3$.
\begin{theorem}

(i) The limiting measure $lim_{\epsilon\to0}\mathcal{D}(I^{\epsilon}(\cdot))$ exists. It is the law of a weak solution 
$I^0(\tau)$ of (\ref{pkdvpr-av}) with the initial data 
$I^0(0)=I(u_0)$.

(ii)The law $\mathcal{D}(I^0(\tau))$ equals to that of $I(v(\tau))$, where $v(\tau)$, $0\leqslant \tau\leqslant T$,  is a  regular solution of the corresponding effective system (\ref{effective-r}) with the  initial data 
$v_0=\Psi(u_0)$. 

(iii) Let $f\geqslant0$ be  a continuous function such that $\int_0^Tf(s)ds=1$. Then the measure $\int_0^Tf(\tau)\mathcal{D}(I^{\epsilon}(\tau),\varphi^{\epsilon}(\tau))d\tau$  on the space $h^p_I\times\mathbb{T}^{\infty}$ 
weakly converges, as $\epsilon \to 0$, to the measure 
$ \Big(\int_0^Tf(\tau)
\mathcal{D}(I^0(\tau))ds \Big)
\times d\varphi  $. 
In particular, the measure 
 $\int_0^Tf(\tau)\mathcal{D}(\varphi^{\epsilon}(\tau))ds$ 
weakly converges to 
$d\varphi$, where $d\varphi$ is the Haar measure on the infinite dimensional torus $\mathbb{T}^{\infty}$.

(iv)  Every sequence $\epsilon^{\prime}_j\to 0$ contains a subsequence $\epsilon_j\to0$ such that the double limit 
$\lim_{\epsilon_j\to0}\lim_{t\to \infty}\mathcal{D}\big(\Psi(u^{\epsilon}(t))\big)
$
exists for any solution $u^{\epsilon}(t)$ and is a stationary measure\footnote{For this notion and its 
discussion see \cite{KS}. We are certain that eq.~(\ref{pkdvp-r}) has a unique stationary measure. When
this is proven, it would imply that the convergence in (iv) holds as $\epsilon\to0$.
}
 for the effective equation (\ref{effective-r}).
\label{thpkdvpr}
\end{theorem}

For the proof see \cite{kpi2008,kuk2010}. For the last assertion also see the argument 
Section 4 in \cite{kuk2011}.  The proof of the theorem applies to other stochastic perturbations of KdV.
In particular, assertions (i)-(iii) hold for equations 
$$
\dot{u}+u_{xxx}-6uu_x=\epsilon g(u(x),x)
+\sqrt{\epsilon}\, \eta,
$$
where $\eta$ is the same as in (\ref{pkdvp-r}) and $g$ is a smooth function, periodic in 
$x$, which has at most a linear growth in $u$, and is such that
$g(u(\cdot), x)\in H^p$ if $u\in H^p$ (this holds e.g. if $g(u,x)$ is even in $u$ and 
odd in $x$).

\medskip
The key to the proof of  Theorem \ref{thpkdvpr}
 is the following result (see Lemma~5.2 in \cite{kpi2008}), where  for any $m\in\mathbb{N}$, $K>0$ and $\delta>0$ we denote 
\begin{equation}\label{Omega}
\fl
\qquad\eqalign{
\Omega(\delta,m,K):=\big\{I: |W_1(I) k_1+& \dots+W_m(I)k_m|<\delta,
\cr
&\mbox{for some}\;k\in\mathbb{Z}^m\;  \mbox{such that} \; 1\leqslant |k|\leqslant K\big\}.
}
\end{equation}

\begin{lem}
For any $m\in\mathbb{N}$, $K>0$, $T>0$ and $\delta>0$ we have 
\[\int_0^T \mathbf{P}\{I(\tau)\in\Omega(\delta,m,K)\}d\tau\leqslant \kappa(\delta,K,m,T),\]
where 
$\kappa(\delta,K, m, T)$  goes to zero with $\delta$,  for any fixed $K$, $m$ and  $T$.
\label{lem-nonresonance}
\end{lem}
This lemma assures that in average,  solutions of (\ref{pkdvp-r}) do not spend much time in the vicinity of   resonant tori.  The stochastic nature of the equation is crucial for  this result. 
\medskip

\subsection{ Deterministic  averaging. }\label{s_determ}
It is  plausible that the averaging principle also holds for equation (\ref{pkdvp1}). But without randomness, it is unclear how to assure that solutions of (\ref{pkdvp1}) `pass the  resonant zone quickly'  (in  analogy with Lemma~\ref{lem-nonresonance}).  This naturally leads to the  question: for which  deterministic perturbations of KdV it is possible to prove
the property of fast crossing the resonant zones and verify the averaging principle?
 Some results in this direction are 
obtained  by the first author in \cite{hg2013, hg20132}. Now we discuss them.

Consider a deterministically perturbed KdV equation: 
\begin{equation}
\dot{u}+u_{xxx}-6uu_x=\epsilon f(u), \quad x\in\mathbb{T}, \;\;u\in H^p,
\label{pkdvdp}
\end{equation}
where $p\geqslant 3$ and the perturbation $f(u)=f(u(\cdot))$ may be non-local. I.e., 
$f(u)(x)$ may depend on values of $u(y)$, where $|y-x|\ge \varkappa>0$.
We  are going to discuss  solutions of (\ref{pkdvdp})
on time-intervals
of order $\epsilon^{-1}$.
Accordingly we fixe some $\zeta_0\le0$, $p\geqslant 3$, $T>0$ and 
make the following assumption:
\smallskip

\noindent{\bf Assumption  A}. { \it(i) There exists  $p'=p^{\prime}(p)<p$, such that for any $q\in [p^{\prime},p]$   the perturbation in (\ref{pkdvdp}) defines an analytic mapping of order $\zeta_0$: 
\[
 H^{q}\to H^{q-\zeta_0},\quad u(\cdot) \mapsto f(u(\cdot)).\]

(ii) For any  $u_0\in H^p$,  there exists a unique solution $u(t)\in H^p$  of (\ref{pkdvdp}) with $u(0)=u_0$.
For $0\le  t \leqslant T\epsilon^{-1}$ its norm  satisfies
$
||u(t)||_p\leqslant C(T,p,||u_0||_p)$.
}

It will be convenient  for us to discuss  equation ({\ref{pkdvdp}})  in the $v$-variables. We will denote
 \[B_p(M)=\{v\in h^p: |v|_p\leqslant M\}.\] With some abuse of notation we will  denote by $S_t$, $0\leqslant t\leqslant T\epsilon^{-1}$, the flow-maps of equation ({\ref{pkdvdp}}), both in the $u$- and in the $v$-variables.

 \begin{defn} 1) A Borelian 
  measure $\mu$ on $h^p$ is called regular if for any analytic function $g\not\equiv0$ on $h^p$, we have $\mu(\{v\in h^p:\;g(v)=0\})=0$.
 
  2) A  measure $\mu$ on $h^p$ is said to be {\it $\epsilon$-quasi-invariant} for equation (\ref{pkdvdp})  if it is regular and for any $M>0$  there exists a constant $C(T,M)$ such that 
  for every Borel set $A\subset B_p(M)$ we have 
\footnote{This specifies the usual definition of a quasi-invariant measure. 
We recall that if a flow $\{S_t\}$ of some equation exists for all $t\ge0$ (or for all $t\in \mathbb R$), 
then a measure $m$ is called quasi-invariant for this equation if the measures $S_t\circ m$ are 
absolutely continuous with respect to $m$ for all $t\ge0$ (respectively for all $t\in\mathbb R$).
}
\begin{equation}
\fl
\quad\quad\quad\quad e^{-\epsilon t C(T,M) }\mu(A)\leqslant \mu(S^t(A))\leqslant e^{\epsilon t C(T,M)}\mu(A),\quad \forall \;
0\le t  \leqslant \epsilon^{-1}T.
\label{quasi-invariant}
\end{equation}
\end{defn}
Similarly, these definitions can be carried to measures on the space $H^p$ and the flow maps of equation (\ref{pkdvdp}) on $H^p$. 

For an $\epsilon$-perturbed finite-dimensional hamiltonian  system the Lebesgue measure is
 $\epsilon$-quasi-invariant 
by the Liouville theorem. This fact is crucial for the Anosov approach to 
  justify the classical averaging principle (see in \cite{AKN, lom1988}). In infinite dimension there is no
 Lebesgue measure,  and existence of
an $\epsilon$-quasi-invariant measure is a serious restriction.

If  equation (\ref{pkdvdp}) has  an $\epsilon$-quasi-invariant measure $\mu$, then the argument, invented by Anosov for the finite dimensional averaging,  insures the required analogy 
of Lemma~\ref{lem-nonresonance} for equation (\ref{pkdvdp}).  Indeed, let us  define  the  resonant
 subset  $  \mathcal{B}$ of $h^p\times\mathbb{R}$ as 
      \[
      \mathcal{B}:=\big\{(v,t): v\in B_p(M),\;\; t\in[0,\epsilon^{-1}T]\;\;\mbox{and}\;\; S^{t}v\in\Omega(\delta,m,K)\big\}
      \]
      (see (\ref{Omega})), and consider the  measure ${\mbox{\boldmath$\mu$}}$ on  $h^p\times \mathbb{R}$,
 where 
 $d {\mbox{\boldmath$\mu$}}=d\mu dt$.
      Then by (\ref{quasi-invariant})  we have 
            \[\fl\quad\quad{{\mbox{\boldmath$\mu$}}}(\mathcal{B})=\int_0^{\epsilon^{-1}T}\mu\Big(B_p(M)\cap S^{-t}\big(\Omega(\delta,m,K)\big)\Big)dt\leqslant \epsilon^{-1}Te^{C(T,M)}\mu(\Omega(\delta,m,K)).\]
      For any $v\in B_p(M)$, define $Res(v)$ as the set of resonant instants of time 
      for a trajectory,  which starts from $v$:
      $$
      Res(v)=\{\tau\in[0,\epsilon^{-1}T]:\;\;S^{t}(v)\in\Omega(\delta,m,K)\}.
      $$
       Then  
      \[{{\mbox{\boldmath$\mu$}}}(\mathcal{B})=\int_{B_p(M)}mes\big(Res(v)\big)d\mu(v)\]
      by the Fubini theorem, 
      where $mes(\cdot)$ stands for the Lebesgue measure on $[0,\epsilon^{-1}T]$. 
   If  for  $\rho>0$ 
 we denote
\[\mathcal{V}\mbox{Res}(\rho):=\{v\in B_p(M): mes\,(Res(v))>\epsilon^{-1}\rho\}\,,
\]
 then in view of the Chebyshev inequality  
\[
\mu(\mathcal{V}\mbox{Res}(\rho))
 \leqslant \frac{\epsilon}{\rho} {\mbox{\boldmath$\mu$}} (\mathcal{B})
\leqslant \frac{Te^{C(M,T)}}{\rho}\mu(\Omega(\delta,m,K)).\]
By Theorem \ref{t_kapp} the functions $v  \mapsto W_1(I(v))k_1+\dots+W_m(I(v))k_m$, $k\in{\mathbb Z}^m\setminus\{0\}$, 
are analytic on $h^p$. Since  they do not vanish identically 
 by Lemma \ref{lem-nd} and the measure $\mu$ is regular, then $\mu(\Omega(0,m,K))=0$.
 Accordingly $\mu(\Omega(\delta,m,K))$  goes to zero with $\delta$, and
 \begin{equation}\label{analogy}
 \mu(\mathcal{V}\mbox{Res}(\rho))\to0 \quad \mbox{as}\quad \delta\to0,
 \end{equation}
 for any $\rho$.   This gives us a sought for  analogy of 
  Lemma~\ref{lem-nonresonance} for deterministic perturbations of KdV which have 
  $\epsilon$-quasi-invariant measures. 
  \medskip

The averaged equation for actions, corresponding to 
 (\ref{pkdvdp}), reads 
\begin{equation}
\frac{dJ_k}{d\tau}=\langle F_k\rangle(J),\quad k=1,2,\dots,
\label{pkdvdp-a}
\end{equation}
where $F_k=\big(d\Psi_k(\Psi^{-1}(v))(f(\Psi^{-1}(v)(\cdot)),\mathbf{v}_k\big)$ (cf. (\ref{pkdvp1-i})).
Due to item (i) of Assumption A, the r.h.s.  of (\ref{pkdvdp-a}) defines a Lipschitz vector filed on $h^p_I$, so the averaged equation is  well posed locally on  $h^p_I$.
We denote by $J_{I_0}(\tau)$ a solution of (\ref{pkdvdp-a}) with an 
initial data $J_{I_0}(0)=I_0\in h^p_I$.
It is shown in \cite{hg2013, hg20132} that 
relation (\ref{analogy}) and the well-posedness of the averaged equation jointly 
allow to establish an averaging theorem  for equation  (\ref{pkdvdp}), provided that it has an
$\epsilon$-quasi-invariant mesure.

 In the statement below 
 $u^\epsilon(t)$ stand for solutions of equation  (\ref{pkdvdp}) and $v^\epsilon(\tau)$ -- for these solutions, 
 written using  the $v$-variables and slow time $\tau=\epsilon t$. 
 By Assumption~A, for $\tau\in[0,T]$ we have 
 $ |I(v^{\epsilon}(\tau))|_p^{\sim}\leqslant C_1\big(T,|I(v^{\epsilon}(0))|^{\sim}_p\big).$ 
Denote 
  $$\tilde{T}(I_0):=\min\{\tau\in\mathbb{R}_+: |J_{I_0}(\tau)|_p^{\sim}\geqslant C_1(T,|I_0|_p^{\sim})+1\}.
  $$
  
\begin{theorem} Fix some $M>0$. 
 Suppose that  Assumption A  \mbox{ holds and equation (\ref{pkdvdp})} has an    $\epsilon$-quasi-invariant measure 
  $\mu$ on $h^p$ such that   $\mu(B_p(M))>0$.  Then
  
  (i) For any $\rho>0$ and  any $q<p-\frac{1}{2}\max\{\zeta_0,-1\}$, there exist $ \epsilon_{\rho,q}>0$ and a Borel subset $\Gamma_{\rho,q}^{\epsilon}\subset B_p(M)$, satisfying
   $\  \lim_{\epsilon\to0}\mu(B_p(M)\setminus \Gamma_{\rho,q}^{\epsilon})=0$,
  with the following property:\\
     if 
   $\epsilon\leqslant \epsilon_{\rho,q}$  and  $ v^{\epsilon}(0)\in \Gamma_{\rho,q}^{\epsilon}$, then
  \begin{equation}
  |I(v^{\epsilon}(\tau))-J_{I_0^\epsilon}(\tau)|_q^{\sim}\leqslant \rho \quad\mbox{for} \quad0\leqslant\tau\leqslant \min\{T,\tilde{T}(I^{\epsilon}_0)\},
 \label{aver}
 \end{equation}
where $I_0^\epsilon=I(v^{\epsilon}(0))$.

(ii) Let $\lambda_{\epsilon}^{v_0}$ be a probability measure on $\mathbb{T}^{\infty}$, defined by the relation 
$$
\fl\quad\quad\int_{\mathbb{T}^{\infty}}f(\varphi)\,d\lambda_{\epsilon}^{v_0}(d\varphi)=
 \frac{1}{T}\int_0^Tf(\varphi(v^{\epsilon}(\tau))d\tau, 
 \quad \forall f\in C(\mathbb{T}^{\infty}),
$$
where $v_0:=v^\epsilon(0)\in B_p(M)$. 
Then the measure 
$ {\mu\big(B_p(M)\big)}^{-1} \int_{B_p(M)}   \lambda_{\epsilon}^{v_0}d\mu(v_0 ) $
  converges weakly, as $\epsilon\to0$,  to the Haar measure $d\varphi$ on $\mathbb{T}^{\infty}$\footnote{In \cite{hg2013} a stronger   assertion  was claimed.
  Namely, that the measure $\lambda_{\epsilon}^{v_0}$ converges to $d\varphi$
  for $\mu$-a.a. $v_0$ in $B_p(M)$. Unfortunately,  the proof  in 
  \cite{hg2013} contains a gap which we still cannot fix. } 
\label{thmpkdvdp}
\end{theorem}

\begin{proposition} \label{prop}
If   Assumption A holds with $\zeta_0<0$,  then for  $\rho<1$, $p\leqslant q$ and $\epsilon\leqslant \epsilon_{\rho,q}$, we have $\tilde T:= \tilde{T}\big(I(v^{\epsilon}(0))\big)>T$ for $v^{\epsilon}(0)\in\Gamma_{\rho,q}^{\epsilon}$.
So  (\ref{aver}) holds for $0\le\tau\le T$. 
\end{proposition}
{\it Proof:}
  Assume that  $\tilde{T} \leqslant T$.  
  By (\ref{aver}) for 
  $0\le\tau\le\tilde T$ 
  we have $|I^{\epsilon}(\tau)-J_{I^{\epsilon}_0}(\tau)|_p^{\sim}\leqslant \rho$.   Therefore $|J_{I_0^{\epsilon}}(\tau )|_p^{\sim}\leqslant C_1(T,|I_0^{\epsilon}|_p^{\sim})+\rho<C_1(T,|I_0^{\epsilon}|_p^{\sim})+1$. This contradicts  the definition of $\tilde{T}$, so $\tilde{T}> T$.\quad $\square$
\smallskip

\noindent{\bf Remark}.  Assume that  $\epsilon$-quasi-invariant measure $\mu$ depends on $\epsilon$, i.e., 
$\mu=\mu_\epsilon $, 
 and  

a)  $\mu_{\epsilon}(\Omega(\delta,m,K))$ goes to zero with $\delta$ \, uniformly in $\epsilon$,

  b)  the constants $C(T,M)$ in (\ref{quasi-invariant}) are bounded  uniformly in $\epsilon$. \\
   Then assertion~(i) holds with $\mu$ replaced by $\mu_{\epsilon}$. For assertion (ii) to hold, more restrictions  should be imposed, see \cite{hg20132}. 
\medskip

Theorem  \ref{thmpkdvdp} gives rise to the questions:

\begin{problem} \label{p.aver}
Does a version of the averaging theorem above holds without assuming  the existence
of an  $\epsilon$-quasi-invariant measure?
\end{problem}

\begin{problem} \label{p.qim}
Which equations   (\ref{pkdvdp}) have $\epsilon$-quasi-invariant measures?
\end{problem}
See the next subsection for some results in this direction.
\begin{problem} \label{higher_order}
Find an averaging theorem 
  for equations  (\ref{pkdvdp}), where the nonlinearity 
defines an unbounded operator, i.e. in  Assumption~A(i) we have $\zeta_0>0$ (note that in the equation from
Example 
in Section~\ref{skam} we have 
$\zeta_0=1$, and in equation (\ref{pkdvp1}) $\zeta_0= 2$). 
\end{problem} 
 It is unlikely that the assertion of Theorem \ref{thmpkdvdp} holds for all initial data, and we believe that the phenomenon of capture  in resonance happens for some solutions of~(\ref{pkdvdp}):
\begin{problem}
Prove that (\ref{aver}) does not hold for some  solutions of (\ref{pkdvdp}).
\end{problem}

\subsection{Existence of $\epsilon$-quasi-invariant measures}
Clearly every regular measure, invariant 
for  equation (\ref{pkdvdp}),  is  $\epsilon$-quasi-invariant. Gibbs measures for some equations of the KdV type are
regular and invariant, they were studied by a number of people (e.g., see \cite{Bo94, Zhi01}).  However,  for generic 
 hamiltonian perturbations of KdV it is difficult, probably impossible,
 to construct  invariant measures in higher order Sobolev spaces due to the lack of high order conservation laws.
Below we give some examples of $\epsilon$-quasi-invariant measures for smoothing 
perturbations of KdV, which are Gibbs measures of KdV (so they are KdV-invariant). Note that some of these 
perturbed equations do not have non-trivial invariant measures. For example, our argument 
applies to equations which in the $v$-variables read as 
$$
\dot{\mathbf{v}}_j=\mathbb{J}W_j(I)\mathbf{v}_j  -\epsilon j^{-\rho}{\mathbf{v}}_j\,,\qquad
\;j\in \mathbb{N},
$$ 
where $\rho>1$. But all trajectories of this equation  converge to zero, so  the only its invariant measure 
 is the $\delta$-measure at the origin.
That is, for averaging in the perturbed KdV (\ref{pkdvdp})
(various) Gibbs measures of KdV play the same role as the Lebesgue measure plays for the classical 
averaging. 

\begin{defn} For any $\zeta_0^{\prime}<-1$, 
a Gaussian measure $\mu_0$ on the Hilbert space $h^p$ 
is called $\zeta^{\prime}_0$-admissible if it has zero mean value
and  a diagonal correlation operator \mbox{$(\mathbf{v}_1,\dots)\mapsto (\sigma_1\mathbf{v}_1,\dots)$}, where
$0<j^{\zeta^{\prime}_0}/\sigma_j\leqslant Const$ for each $j$. 
\end{defn}

For  any $\zeta_0^{\prime}<-1$  a $\zeta_0^{\prime}$-admissible measure   $\mu_0$ is a well-defined probability measure on $h^p$, which can be formally written 
as 
 \begin{equation}
 \mu_0=\prod_{j=1}^{\infty}\frac{(2\pi j)^{1+2p}}{2\pi\sigma_j}\exp\{-\frac{(2\pi j)^{1+2p}|\mathbf{v}_j|^2}{2\sigma_j}\}d\mathbf{v}_j,
 \label{gaussian}
 \end{equation}
  where $d\mathbf{v}_j$, $j\geqslant 1$, is the Lebesgue measure on 
  $\mathbb{R}^2_{\mathbf{v}_j }$. It is known  that (\ref{gaussian})   is a well-defined  measure on $h^p$ if and only if  $\sum \sigma_j<\infty$ (see \cite{Bogachev}). It  is regular and non-degenerate 
  in the sense that its support equals $h^p$ (see \cite{Bogachev, bom2013}). Writing KdV in the $v$-variables   
    we see that $\mu_0$ is  invariant under the KdV flow.  Note that $\mu_0$ is a Gibbs measure for KdV, written in the
    form (\ref{bnf-e}), since formally it may be written as $\mu_0=Z^{-1}\exp\{-\langle Qv,v\rangle\}dv$, where $\langle Q v,v\rangle=\sum c_j|\mathbf{v}_j|^2$ is an integrals of motion for KdV (the statistical sum $Z=\infty$, so indeed this is a formal expression).
  
  For a perturbed KdV (\ref{pkdvdp}) we define $\mathcal{P}(v)=d\Psi(u)(f(u))$, where $u=\Psi^{-1}(v)$. 
  A non-complicated 
  calculation (see in \cite{hg2013}) shows that:
\begin{theorem}If Assumption A holds and \\
$(i)^{\prime}$ the operator $\mathcal{P}$ analytically maps the space $h^p$ to $h^{p-\zeta_0^{\prime}}$ with some $\zeta_0^{\prime}<-1$, \\
then  every $\zeta_0^{\prime}$-admissible Gaussian measure on $h^p$ is 
$\epsilon$-quasi-invariant for equation~(\ref{pkdvdp}) on the space $h^p$.
\label{quasi-measure1}
\end{theorem}

However, due to the complexity of the nonlinear Fourier transform $\Psi$, it is not easy to verify the condition $(i)^{\prime}$ of Theorem \ref{quasi-measure1} for specific equation (\ref{pkdvdp}). Now we will give other examples of $\epsilon$-quasi-invariant measures on the space $H^p$, by  strengthening the  restrictions in Assumption~A.
Suppose that there $ p\in\mathbb{N}$. Let $\mu_p$ be the centered Gaussian measure on $H^p$ with the 
correlation operator $\triangle^{-1}$. Since $\triangle^{-1}$ is an operator of the trace type, then $\mu_p$ is a well-defined probability measure on $H^p$. 

We recall (see Remark \ref{remark-laws}) that KdV has infinitely  many conservation laws $\mathcal{J}_n(u)$, $n\geqslant0$, of the form $\mathcal{J}_n=\frac{1}{2}||u||_n^2+J_{n-1}(u)$, where $J_{-1}(u)=0$ and for $n\geqslant1$,
\begin{equation}
\fl\quad\quad\quad\quad
J_{n-1}(u)=\int_{\mathbb{T}}\big\{c_nu(\partial_x^{n-1}u)^2+\mathcal{Q}_n(u,\dots,\partial^{n-2}_xu)\big\}dx\,.
\label{law-form}
\end{equation}
 Here  $c_n$ are real constants and $\mathcal{Q}_n$ are polynomial in their arguments.
  From (\ref{law-form}), we know that  the functional $J_p$ is  bounded on bounded sets in $H^p$.
We consider a Gibbs measure $\eta_p$ for KdV,  defined by its density against~$\mu_p$, $$
\eta_p(du)=e^{-J_p(u)}\mu_p(du).$$
It  is  regular and non-degenerated in the sense that its support contains the whole space $H^p$ (see  \cite{Bogachev}). Moreover, it is invariant for KdV \cite{Zhi01}. The following theorem was shown in \cite{hg20132}:
\begin{theorem}
If Assumption A holds with $\zeta_0\leqslant -2$, then the measure $\eta_p$ is $\epsilon$-quasi-invariant for perturbed KdV (\ref{pkdvdp}) on the space $H^p$.
\end{theorem}

\begin{corollary} If $\zeta_0\leqslant -2$, then the  assertions of Theorem~\ref{thmpkdvdp} hold with $\mu=\eta_p$, and we have 
$\tilde T>T$. 
\label{coro-th}
\end{corollary}

In particular, this corollary applies to the equation
\begin{equation*}
\dot{u}+u_{xxx}-6uu_x=\epsilon f(x), \quad x\in\mathbb{T}, \;\;u\in H^p,
\end{equation*}
where $f(x)$ is a smooth function with zero mean-value.  
This equation  may be viewed as a model for shallow water wave propagation under small external force. Note that the KAM-Theorem~\ref{kdv-kam} also applies to it.

\begin{problem}
Besides the class of  $\zeta^{\prime}_0$-admissible
 Gaussian measures and Gibbs measure $\eta_p$,  there are   many other KdV-invariant measures. How to check 
  if a measure like that is   $\epsilon$-quasi-invariant for a given $\epsilon$-perturbation of KdV?
\end{problem}

\subsection{  Nekhoroshev type results (long-time stability)?}
 In  the finite dimensional case, the strict convexity of the unperturbed integrable Hamiltonian assures the long-time stability of  solutions for  perturbed hamiltonian  equations (\cite{nek1972, loc1992, jpo1993}).  Theorem \ref{th-convex} tells us the the KdV Hamiltonian $H_K(I)$ is convex in $l^2$ and hints that it is strictly convex (at least) in a neighborhood of the origin in $l^2$.  This suggests that an Nekhoroshev type stability may hold for perturbed KdV under hamiltonian  perturbations  (see equation (\ref{pkdvh1})), at least for initial data in a neighborhood of the origin, where the strict convexity should hold. But at the moment
 of writing   no exact statement is available.
 
 There are  several {\it ad hoc} quasi-Nekhoroshev theorems for hamiltonian  PDEs, see 
 \cite{ bam1999-2,  bou2000} and references therein.  However, these results only apply in a small  (of the size of the perturbation) neighborhood of the origin. Nonetheless, we believe that the corresponding  technique and the results in Theorem \ref{th-convex} will lead to results on long-time stability (at least in time interval of order $\epsilon^{-p}, p\ge1$,) for some solutions of perturbed KdV (\ref{pkdvh1}). Note that 
  Theorem~\ref{thmpkdvdp} and Corollary \ref{coro-th} imply  such a stability for $p=1$ and for typical initial data, if the perturbation 
  $f$ is 
  smoothing. 
  Stability on time-intervals of order  $\epsilon^{-2}$ seems to be a 
  much harder question. 
 \begin{problem}
 Is it possible to prove a long-time stability result for perturbed KdV under hamiltonian  perturbations, e.g. for equation (\ref{pkdvh1}), that holds for all `smooth' initial data? 
 \end{problem}

\noindent 
  {\bf Acknowledgments.} \small This work was supported by
 l'Agence Nacionale de la Recherche, grant ANR-10-BLAN 0102.

\section*{References}
 \bibliography{review_kdv.bib}

\end{document}